%% file: paper.tex
\documentclass[11pt]{article}

\usepackage[utf8]{inputenc}
\usepackage[T1]{fontenc}

\usepackage{amsmath,amssymb,amsfonts,textcomp,amsthm,xifthen,graphicx,color,pgfplots}

\usepackage{enumerate}
\usepackage{enumitem}

\usepackage{fullpage}

\usepackage[utf8]{inputenc}

\usepackage[pdftex,
            pdfauthor={Thomas F\"uhrer and Norbert Heuer},
            pdftitle={Mixed finite elements for Kirchhoff--Love plate bending},
            ]{hyperref}

\usepackage{pgfplots}
\usepgfplotslibrary{groupplots}
\pgfplotsset{compat=newest}
\usepgfplotslibrary{external}
\usepgfplotslibrary{colorbrewer}
\input{header}

\title{Mixed finite elements for Kirchhoff--Love plate bending
\thanks{Supported by ANID through FONDECYT projects 1210391, 1230013}}
\author{
Thomas~F\"uhrer$^\dagger$
\and
Norbert Heuer\thanks{
Facultad de Matem\'aticas, Pontificia Universidad Cat\'olica de Chile,
Avenida Vicu\~na Mackenna 4860, Santiago, Chile,
email: {\tt \{tofuhrer,nheuer\}@mat.uc.cl}}}

\date{}
\begin{document}
\maketitle

\begin{abstract}
We present a mixed finite element method with triangular and parallelogram meshes for the Kirchhoff--Love
plate bending model. Critical ingredient is the construction of low-dimensional local
spaces and appropriate degrees of freedom that provide conformity
in terms of a sufficiently large tensor space and that allow for
any kind of physically relevant Dirichlet and Neumann boundary conditions.
For Dirichlet boundary conditions and polygonal plates, we prove quasi-optimal convergence of the
mixed scheme.
An a posteriori error estimator is derived for the special case of the biharmonic problem.
Numerical results for regular and singular examples illustrate our findings.
They confirm expected convergence rates and exemplify the performance of an
adaptive algorithm steered by our error estimator.

\end{abstract}
%===================================================================================================

%===================================================================================================
\section{Introduction}

Plate bending models have been the subject of research in numerical analysis for several decades,
until today. This is not only due to their relevance in structural engineering
but also owed to the inherent mathematical challenges.
The Kirchhoff--Love and Reissner--Mindlin models are the classical ones.
The former can be interpreted as the singularly perturbed
limit of the latter for plate thickness tending to zero. 
This limit case poses the vertical deflection $u$ as an $H^2(\Omega)$-function whereas
the bending moments $\bM$ are set in the space of symmetric $L_2(\Omega)$-tensors with
$\div\Div\bM\in L_2(\Omega)$ for $L_2$-regular vertical loads $f$, in short,
$\bM\in\HdivDivset{\Omega}$.
(Here, $\Omega\subset\R^2$ denotes the plate's mid-surface, and
$\Div\bM$ means the row-wise application of the divergence operator.)
For non-convex polygonal plates, $\HdivDivset{\Omega}$ is not a subspace of $\HH^1(\Omega)$
(tensors with $H^1(\Omega)$-components) or $\HDivset{\Omega}$ (symmetric tensors $\bM$ with
$\Div\bM\in \LL_2(\Omega):=L_2(\Omega)^2$), cf.~\cite{BlumRannacher80}.
This lack of regularity constitutes a serious challenge
for the approximation of bending moments and its analysis.
Bending moments are critical quantities in engineering applications and have been elusive
to conforming approximations in $\HdivDivset{\Omega}$ for a long time. An early, only partially
conforming approach is the Hellan--Herrmann--Johnson method that gives bending
moment approximations with continuous normal-normal traces,
see~\cite{Hellan_67_AEP,Herrmann_67_FEB,Johnson_73_CMF}.

In this paper, we present low-dimensional finite elements of low degree
on triangles and parallelograms for bending moments with corresponding basis functions, and
prove quasi-optimal convergence of the corresponding mixed finite element scheme.
This includes the critical case of non-convex polygonal plates.
For triangles $K$, the space for the bending moments is
\begin{align*}
  \XXX(K) = \symCurl(\RT^2(K))\oplus \bx\bx^\top P^1(K) = \sym(\RT^0(K)\otimes\RT^1(K))
\end{align*}
and has dimension $15$.
The construction for parallelograms is similar and leads to a space of dimension $20$.
Here, $P^p(K)$ are polynomials of degree $\leq p$, $\RT^p(K) = \bP^p(K) + \bx P^p(K)$
is the Raviart--Thomas space of order $p$ and
$\symCurl(\cdot) = \tfrac12(\Curl(\cdot)+\Curl(\cdot)^\top)$.

We define degrees of freedom that provide $\HdivDivset{\Omega}$-conformity without requiring any
additional non-physical regularity. Let us discuss this in some detail.
Here and in previous papers, our focus is on discretizations of plate bending moments
that a) are \emph{conforming} in $\HdivDivset{\Omega}$,
b) have \emph{few degrees of freedom},
c) are \emph{explicit} with respect to a basis,
d) allow for the \emph{least regular cases} appearing in engineering applications,
and e) provide \emph{optimal approximation orders} under these conditions.
Indeed, regularity is a critical issue in plate bending analysis.
The Kirchhoff--Love model (considered here) allows for bending moments with
jumps at incoming corners. These \emph{corner forces}
can be characterized as delta distributions acting on deflections and
appear in engineering applications, causing notorious problems in computational mechanics.

The $\HdivDivset{\Omega}$-element of this paper is the first that combines all the
mentioned requirements in a finite element setting, and is based on our previous studies
of the DPG method \cite{KLove1,KLove2}.
Specifically, in \cite{KLove1} we gave the first well-posed local (edge-wise)
characterization of conformity in a subspace of $\HdivDivset{\Omega}$
that includes bending moments representing corner forces,
and identified the smallest set of degrees of freedom that provide conformity for a
discrete setting based on traces.
Those findings comply with all the requirements above.
(The approximation order there is one but can be increased to two by adding
one degree of freedom per edge.)
Here, we develop and analyze bending moment elements for triangles and
parallelograms. The degrees of freedom coincide with those of our DPG trace
approximation (enriched to second order).
Our elements comply with all the requirements a)--d) and provide an
$\HdivDivset{\Omega}$ interpolation operator that is well defined on a
subspace that includes corner forces.
In particular, we provide an explicit basis that is dual to the degrees
of freedom, including a basis-function representation of corner forces.
Having a basis for the degrees of freedom, an implementation is relatively
straightforward. In fact, to our knowledge, the implementations of DPG
from \cite{KLove1}, with extensions reported in \cite{KLove2,FuehrerHH_21_TOB,FuehrerHN_22_DMS},
and of the mixed finite element method from this paper
are so far the only conforming $\HdivDivset{\Omega}$ implementations that include
corner forces.
Let us stress the fact that our characterization of conformity and the corresponding
degrees of freedom coincide with natural boundary data for Kirchhoff--Love
plate bending:
normal-normal bending moments, effective shear forces, and corner forces.

Independently of plate bending models,
there is wider interest in approximations of $\HdivDivset{\Omega}$ for
domains $\Omega\subset\R^d$, $d\ge 2$, usually motivated by the bi-Laplacian.
There is a series of papers mainly by Chen and Huang, and some others, on this subject.
In contrast to our low-order 2d-focus, Chen and Huang develop higher-order approximations
in two and more space dimensions.
The two-dimensional case of triangles is studied in
\cite{ChenH_FED}, see also \cite{HuMZ_21_FMF},
and generalized in \cite{ChenHuang22Siam} to simplices in arbitrary space dimensions.
These papers consider degrees of freedom that include point values of bending moments.
They are non-physical in the sense that, even though being of
interest in applications (for sufficiently smooth solutions), they are not
available as data (there is no point-wise physical definition) and
they do not give rise to bounded functionals on
a sufficiently large subspace of $\HdivDivset{\Omega}$ that includes bending
moments of corner forces.
In their 3d-paper \cite{ChenH_22_FEDb}, the authors switch to our formulation
of $\HdivDivset{\Omega}$-conformity from \cite{KLove1} (which is valid in two
and three dimensions), but again consider point values of bending moments.
More recently, in \cite{ChenHuang23}, Chen and Huang take the degrees
of freedom from \cite{KLove1} (thus avoiding point values) and extend them
to three dimensions and higher polynomial degrees.
They propose to use the space $\mathbb{Y}(K) = \PPPsym^k(K)$ with degree $k\geq 3$,
giving dimension $30$ for the lowest-order case $k=3$.
An analogy for the difference between $\XXX(K)$ and $\mathbb{Y}(K)$ is the difference between
Raviart--Thomas and Brezzi--Douglas--Marini elements for discretizations of $\Hdivset\Omega$-vector fields.
Both spaces $\XXX(K)$ and $\mathbb{Y}(K)$ have approximation order two in $\HdivDivset{K}$.
In~\cite[Section~2.4]{ChenHuang23}, a reduced lowest-order element with $20$ degrees of freedom
is presented. However, its implementation requires a basis that is dual to the degrees of freedom
``which can be a challenging task'', to quote from~\cite[p.14]{ChenHuang23}.
In their new preprint \cite{ChenHuang23b}, which is a revision of
\cite{ChenHuang23} with new title, Chen and Huang use the construction of finite element spaces
from this paper, and extend it to higher dimensions exclusively for simplices,
see \cite[\S 2.5]{ChenHuang23b}.
Rather than constructing a basis, Chen and Huang propose a hybridization technique
to avoid a direct mixed FEM implementation, and no numerical results are reported.

We mention that there are some classical mixed schemes for Kirchhoff--Love plates, see, e.g.,
\cite{Falk_78_ABE,BrambleF_83_TMF,Monk_87_MFE}.
They are based on the interpretation of the Kirchhoff--Love model
(with isotropic homogeneous material) as the bi-Laplacian, $\Delta^2 u=f$, and introducing
$v=\Delta u$ as an independent variable, as proposed by Ciarlet and Raviart in 
\cite{CiarletR_74_MFE}. This strategy requires more regularity than generally available
($v\in H^1(\Omega)$) and does not allow for general Neumann boundary conditions as $v$ is
a non-physical variable. We remark that in \cite{FuehrerHH_21_TOB},
we presented a DPG-setting for the two-variable setting with $v=\Delta u$ that is well posed
for non-convex domains and $L_2(\Omega)$-loads.

Our analysis employs techniques and tools that we have learned from our studies \cite{KLove1,KLove2}
of the discontinuous Petrov--Galerkin (DPG) method in the context of plate bending.
Whereas the DPG framework may seem to be very specialized and irrelevant for the analysis
of classical Galerkin approaches including mixed schemes, we here illustrate that this view
is not correct. The most common DPG setting is based on ultraweak formulations. Their analysis
requires a specific formulation of trace operators and the discretization of their images,
the resulting trace spaces. On the domain level, these traces give rise to precise
conditions of conformity, e.g., across interfaces. For canonical spaces this is well known.
For example, $H^1(\Omega)$-conformity requires continuity in the sense of $H^{1/2}$-traces and
$\Hdivset{\Omega}$-conformity means the $H^{-1/2}$-continuity of normal components.
There are spaces where such an approach to conformity is much more intricate,
for instance $\HdivDivset{\Omega}$ introduced before.
We stress the fact that there is a key difference between the conformity in the full space and
the conformity of piecewise polynomial (or otherwise) approximations.
In the latter case, trace operators have to be localized. Considering that trace spaces
are typically of fractional order, this is a serious challenge. To circumvent this problem
one usually requires more regularity. For instance, normal traces of $\Hdivset{\Omega}$
are considered in $L_2$ rather than $H^{-1/2}$. The key point is to increase the regularity
as little as possible in order not to exclude relevant cases of low regularity.
In this paper, we present and analyze a conforming piecewise polynomial approximation
of bending moments $\bM$ which only requires a slightly increased regularity
$\bM\in\HdivDivset{\Omega,\mathcal{E}}$.
Here, $\HdivDivset{\Omega,\mathcal{E}}\subset \HdivDivset{\Omega}$ is a dense
subspace (introduced in \S\ref{sec:notations}) and therefore, our construction is applicable without
any additional regularity requirement.

We use the trace formulation from \cite{KLove1} to construct
$\HdivDivset{\Omega}$-conforming elements. In contrast to our DPG-setting,
where we use lowest-order moments for both the normal-normal traces and the effective shear forces
(plus vertex jumps of tangential-normal traces), here we additionally use
their first-order moments.
In this way, second-order approximations are achieved
(though, as already mentioned, second order can be easily achieved also in the DPG setting).
There is an inherent piecewise polynomial $\HdivDivset{\Omega}$-interpolation operator.
It commutes with the $L_2$-projection onto piecewise linear polynomials.
Therefore, we have the canonical ingredients to set up a mixed finite element scheme
and prove second order convergence for sufficiently regular solutions.
The vertical deflection is approximated in $L_2(\Omega)$ by piecewise linear polynomials
and the bending moments are approximated in $\HdivDivset{\Omega}$ by our basis functions.
A postprocessing scheme with piecewise cubic approximations of the vertical deflection
is studied as well.

We also derive an a posteriori error estimator for our mixed scheme.
For simplicity we consider constant material properties which is equivalent to
studying the biharmonic problem. An extension to piecewise constant coefficients is possible
but not pursued here.
Error estimators are critical for steering adaptive algorithms in order to resolve singularities. 
In the analysis we follow ideas and techniques from Carstensen~\cite{CC97} for Poisson-type problems.
The main ingredient for proving reliability of the estimator is a Helmholtz decomposition
for vector fields. Here, we consider a Helmholtz-type decomposition of symmetric tensors fields.
It uses an $\symCurl(\cdot)$-representation of tensors which are in the kernel of the $\dDiv$ operator,
see, e.g.,~\cite{RafetsederZulehner18} and references therein.
We stress the fact that our error estimator is not restricted to using space $\XXX(K)$
but can be derived for other conforming discretizations as mentioned before. 

An overview of the remainder is as follows.
In Section~\ref{sec:notations} we introduce Sobolev spaces and norms, recall plate-specific
trace operators and their properties, and discuss polynomial spaces and transformation tools.
Fundamental to this paper is our conformity characterization by Proposition~\ref{prop_conf}.
Local and global finite element spaces are the subject of Section~\ref{sec:femSpace}.
Outcome is a piecewise polynomial interpolation operator that provides
$\HdivDivset{\Omega}$-conforming, second-order approximations (Proposition~\ref{prop:proj}).
The mixed formulation and finite element scheme for the Kirchhoff--Love
plate being model are presented in Section~\ref{sec:mixed}. Theorem~\ref{thm:mixed} establishes
the quasi-optimal convergence of the scheme, and in Section~\ref{sec:postproc} we
present and analyze a postprocessing scheme for the primal variable.
Theoretical results in Section~\ref{sec:mixed} are proved for Dirichlet boundary conditions.
Though, we stress the fact that our setting allows for implementing Neumann and
mixed boundary conditions as well. Corresponding proofs require more technical details and
are left open here.
In Section~\ref{sec:aposteriori} we present an a posteriori error estimator
for the bending moments, the Hessian of the deflection in the biharmonic case,
and prove its reliability and local efficiency.
Finally, in Section~\ref{sec:ex} we report on numerical experiments.
They include the case of a non-convex domain and singular solution, and underline the
performance of an adaptive scheme that is based on our error estimator.
In Appendix~\ref{sec:appendix}, we illustrate the construction of basis (shape) functions
for triangles and parallelograms. This is useful for implementation.

\section{Preparations}\label{sec:notations}
\subsection{Sobolev spaces and trace operators}
For a Lipschitz domain $\omega\subset \R^2$ we denote by $L_p(\omega)$ ($p\in[1,\infty]$)
resp. $H^s(\omega)$ ($s\in[0,2]$) Lebesgue resp. Sobolev spaces.
The Sobolev spaces are defined by real interpolation between $L_2(\omega)$ and $H^2(\omega)$.
The canonical norm and inner product in $L_2(\omega)$ are denoted by $\norm{\cdot}{\omega}$
and $\ip{\cdot}\cdot_\omega$. For spaces of vector-valued resp. tensor-valued functions
we use boldfaced resp. blackboard boldfaced letters, e.g., $\LL_2(\omega)=L_2(\omega;\R^2)$
and $\mathbb{L}_2(\omega)=L_2(\omega;\R^{2\times 2})$.
For spaces of symmetric tensor-valued functions we use $\LLLts(\omega)$.
For a boundary part $\gamma\subseteq \partial \omega$ we use $\dual{\cdot}{\cdot}_\gamma$
to denote either a duality or the $L_2(\gamma)$ inner product. 

We follow~\cite{KLove1} and introduce
\begin{align*}
  \HdivDivset\omega := \set{\bM\in\LLLts(\omega)}{\dDiv\bM\in L_2(\omega)},
\end{align*}
a Hilbert space with (squared) norm
\begin{align*}
  \norm{\bM}{\HdivDivset\omega}^2 = \norm{\bM}\omega^2 + \norm{\dDiv\bM}\omega^2.
\end{align*}
Here, $\Div$ is the divergence operator applied row-wise.
We also need $\Curl\bq$ which is the $\curl$ operator applied to each component of $\bq$, and
\begin{align*}
  \symCurl\bq = \frac12\big(\Curl\bq + (\Curl\bq)^\top\big)
\end{align*}
is the symmetrized $\curl$ operator. 
In addition, $\Rot$ is the row-wise $\rot$ operator.
Sobolev space $H^2(\omega)$ is equipped with the (squared) norm
\begin{align*}
  \norm{z}{H^2(\omega)}^2 = \norm{z}\omega^2 + \norm{\Ggrad z}\omega^2.
\end{align*}
We introduce trace operators $\trdDiv{\omega}\colon \HdivDivset\omega\to H^2(\omega)'$ and
$\trGgrad{\omega}\colon H^2(\omega) \to \HdivDivset\omega'$ by
\begin{subequations}\label{eq:def:traceop}
\begin{align}
  \trdDiv\omega(\bM)(z) &:= \dual{\trdDiv\omega\bM}z_{\partial\omega} := \ip{\dDiv\bM}{z}_{\omega} -\ip{\bM}{\Ggrad z}_{\omega}, \\
  \trGgrad\omega(z)(\bM) &:= \dual{\trGgrad\omega z}{\bM}_{\partial\omega} :=  \dual{\trdDiv\omega\bM}z_{\partial\omega}
\end{align}
for all $\bM\in\HdivDivset\omega$, $z\in H^2(\omega)$.
\end{subequations}
While these operators are defined via volume integrals it can be easily seen with integration by parts
that they reduce to boundary terms for sufficiently smooth arguments.

For a Lipschitz domain $\Omega$ let $\TT$ denote a regular decomposition into open triangles and
parallelograms $K$ such that
\begin{align*}
  \overline\Omega = \bigcup_{K\in\TT} \overline K.
\end{align*}
Here, by regular we mean that all elements are non-degenerate and two distinct but touching elements either share one vertex or one edge. In particular, there are no hanging nodes.
The set of edges of the mesh is denoted by $\edges$ and the set of vertices by $\vertices$. We use $\edges_K$, $\vertices_K$ to denote the edges, vertices of an element $K\in\TT$. Furthermore, $\vertices_\Omega$ means the set of all interior vertices and $\vertices_\Gamma$ are the boundary vertices. The analogous notation is used for edges.
The set $\patch_\TT(z) = \set{K\in\TT}{z\in\vertices_K}$ is the neighborhood of vertex $z$ with corresponding domain $\Patch_\TT(z)$. For elements $K\in\TT$ we use the notation $\patch_\TT(K)=\bigcup_{z\in\vertices_K} \patch_\TT(z)$ and $\Patch_\TT(K)$ for the corresponding domain.
The triangular reference element $\Kref_\triangle$ is the interior of the convex hull of vertices
$(0,0)$, $(1,0)$, and $(0,1)$, and the square reference element $\Kref_\square$ is the interior
of the convex hull of vertices $(0,0)$, $(1,0)$, $(1,1)$, and $(0,1)$.
In what follows we use the symbol $\Kref$ to denote either $\Kref_\triangle$ or $\Kref_\square$.
The diameter of an element $K$ resp. edge $E$ is denoted by $h_K$ resp. $h_E$. We define the mesh-size function $h_\TT\in L_\infty(\Omega)$ by $h_\TT|_K = h_K$. Moreover, we assume that $\TT$ is a decomposition of $\Omega$ into shape-regular elements. This implies that neighboring elements have comparable diameters as well as $h_K \eqsim h_E$ for all $E\in\edges_K$ and all $K\in\TT$. 

Here, and for the remainder of this work $A\eqsim B$ means that $A\lesssim B$ and $B\lesssim A$. The estimate $A\lesssim B$ for $A,B\geq 0$ is an abbreviation of $A\leq C\cdot B$ where $C>0$ is a generic constant possibly depending on the shapes of the elements in $\TT$ but independent of their diameter.

For an element $K\in\TT$ we use the generic notation $\tangential$ resp. $\normal$ for the tangential resp. normal vector along $\partial K$ (in positive orientation).
For sufficiently smooth symmetric tensors $\bM$ defined on $K\in\TT$, we set for $E\in\edges_K$
\begin{align*}
  \trdDiv{K,E,\normal}\bM &= \normal\cdot\bM\normal|_E, \\
  \trdDiv{K,E,\tangential}\bM &= (\normal\cdot\Div\bM + \partial_{\tangential}(\tangential\cdot\bM\normal))|_E.
\end{align*}
By $\partial_\tangential$ resp. $\partial_\normal$ we denote tangential resp. normal derivatives along an edge.
For $E,E'\in \edges_K$, $E\neq E'$, $\{z\} = \overline{E}\cap\overline{E'}$ with $z$ being the endpoint of $E$ and starting point of $E'$ we set
\begin{align*}
  \jump{\tangential\cdot\bM\normal}_{\partial K}(z) &= (\tangential\cdot\bM\normal)|_{E}(z) - (\tangential\cdot\bM\normal)|_{E'}(z).
\end{align*}
The latter trace terms can be interpreted as bounded functionals. To that end consider the spaces
\begin{align*}
  H_K^{3/2}(E) := \set{z|_E}{z\in H^2(K)}, \quad H_K^{1/2}(E) := \set{\partial_\normal z|_E}{z\in H^2(K)} \quad\forall E\in\edges_K,
\end{align*}
and
\begin{align*}
  &\HdivDivset{\TT,\edges} 
  \\&\quad:= \set{\bM}{\bM|_K\in\HdivDivset{K}, \trdDiv{K,E,\normal}\bM\in H_K^{1/2}(E)', \, \trdDiv{K,E,\tangential}\bM\in H_K^{3/2}(E)', \,E\in\edges_K,\,K\in\TT}.
\end{align*}
The jump term can then be interpreted as the functional
\begin{align} \label{eq:jump}
  \dual{\jump{\tangential\cdot\bM\normal}_{\partial K}}{v}_{\partial K} &:=
  \sum_{E\in\edges_K} (\dual{\trdDiv{K,E,\tangential}\bM}{v}_E - \dual{\trdDiv{K,E,\normal}\bM}{\partial_\normal v}_E)-\dual{\trdDiv{K}\bM}v_{\partial K}.
\end{align}
For $v\in H^2(K)$ with $v(z) = 1$ and $v(z') =0$ for $z,z'\in \vertices_K$, $z\neq z'$, this gives
\begin{align*}
  \jump{\tangential\cdot\bM\normal}_{\partial K}(z) = \jump{\tangential\cdot\bM\normal}_{\partial K}(z)v(z) := \dual{\jump{\tangential\cdot\bM\normal}_{\partial K}}{v}_{\partial K}.
\end{align*}
The next result characterizes conformity of $\HdivDivset\Omega$ elements, see~\cite[Proposition~3.6]{KLove1}.

\begin{proposition} \label{prop_conf}
  Let $\bM\in \HdivDivset{\TT,\edges}$. Then, $\bM\in\HdivDivset\Omega$ if and only if for all $E\in\edges_\Omega$, $z\in \vertices_\Omega$
  \begin{align*}
    \trdDiv{K,E,\normal}\bM + \trdDiv{K',E,\normal}\bM &= 0, \quad  \trdDiv{K,E,\tangential}\bM + \trdDiv{K',E,\tangential}\bM = 0, \\
    \sum_{K\in\patch_\TT(z)} \jump{\tangential\cdot\bM\normal}_{\partial K}(z) &= 0
  \end{align*}
  where $K,K'\in\TT$, $K\neq K'$ with $\overline K\cap \overline{K'} = \overline E$.
\end{proposition}

\subsection{Polynomial spaces and basis functions}
Polynomial spaces are denoted by $P^p(K)$ (polynomials of degree $\leq p$ on element $K$), $P^p(E)$ (polynomials of degree $\leq p$ on edge $E\in \edges_K$). Vector-valued polynomials are denoted by $\bP^p(K)$, and tensor-valued polynomials by $\PPP^p(K)$. Furthermore, the space of symmetric tensor-valued polynomials is $\PPPsym^p(K)$.
The $L_2(\Omega)$-orthogonal projection onto $P^p(K)$ is denoted by $\Pi_K^p$,
and $\Pi_\TT^p$ is the $L_2(\Omega)$-orthogonal projection onto $P^p(\TT)$.

The construction of our element is based on the Raviart--Thomas space
\begin{align*}
  \RT^p(K) = \bx P_\mathrm{hom}^p(K) \oplus \bP^p(K).
\end{align*}
Here, $P_\mathrm{hom}^p(K)$ denotes the space of homogeneous polynomials of degree $p$. 
We slightly abuse notation by writing $\bx$ for the function $(x,y)\mapsto(x,y)^\top$ while $\bx P_\mathrm{hom}^p(K) = \set{\bx\eta}{\eta\in P_\mathrm{hom}^p(K)}$ denotes a space.
Similarly, $\bx^\top$ corresponds to the function $(x,y)\mapsto (x,y)$.
We stress that $\RT^p(K)$ is well defined for $K$ being a triangle or a parallelogram.
However, for a parallelogram $K$, $\RT^p(K)$ is not the Raviart--Thomas space on parallelograms. We introduce the lowest-order Raviart--Thomas space $\bQ^0(K)$ for parallelograms below in Section~\ref{sec:localSpace:rect}.

For $K\in\TT$ and $E\in\edges_K$, $\ell_{E,k}\in P^k(E)$ denotes the Legendre polynomial of degree $k\in\N_0$ that is normalized so that $\ell_{E,k}(z) = 1$ for $z$ the end point of $E$.

The dyadic product of two column vector-valued functions with the same domain is
\begin{align*}
  \bphi\bpsi^\top\colon z\mapsto\bphi(z)\bpsi(z)^\top.
\end{align*}
Let $\bX$, $\bY$ denote spaces (of vector-valued functions over the same domain), then
\begin{align*}
  \bX\otimes \bY = \linhull\set{\bphi\bpsi^\top}{\bphi\in\bX, \, \bpsi\in\bY}.
\end{align*}
We introduce the symmetrize operation for tensor-valued functions,
\begin{align*}
  \sym(\bM) = \frac12(\bM+\bM^\top)
\end{align*}
and also adopt the notation for spaces, i.e., 
\begin{align*}
  \sym(\XXX) = \set{\sym(\bM)}{\bM\in\XXX}.
\end{align*}
We also need the space
\begin{align*}
  \widetilde\PPP^2(K) = \symCurl(P_\mathrm{hom}^2(K)\bx)
\end{align*}
and note that $\widetilde\PPP^2(K)\subset\ker(\dDiv)$. Furthermore, one verifies that 
$\widetilde\PPP^2(K)$ is spanned by the three tensor-valued functions
\begin{align}\label{eq:tildeP2span}
  (x,y)\mapsto\begin{pmatrix}
    x^2 & 0 \\ 0 & -y^2
  \end{pmatrix}, \quad (x,y)\mapsto
  \begin{pmatrix} xy & \tfrac{y^2}2 \\ \tfrac{y^2}2 & 0
  \end{pmatrix}, \quad (x,y)\mapsto
  \begin{pmatrix} 0 & \tfrac{x^2}2 \\ \tfrac{x^2}2 & xy
  \end{pmatrix}.
\end{align}
Finally, we introduce $\widetilde\PPP^3(\Kref)$ as the space spanned by
\begin{align*}
  (x,y)&\mapsto \begin{pmatrix} x^2 & xy \\ xy & -2y^2 \end{pmatrix}, 
  \quad
  (x,y)\mapsto \begin{pmatrix} 0 & y^2 \\ y^2 & 0 \end{pmatrix}, 
  \quad
  (x,y)\mapsto \begin{pmatrix} 0 & 0 \\ 0 & xy \end{pmatrix}, \\
  \quad
  (x,y)&\mapsto \begin{pmatrix} 2x^3 & -x^2y \\ -x^2y & -4xy^2 \end{pmatrix}, 
  \quad
  (x,y)\mapsto \begin{pmatrix} 4x^2y & xy^2 \\ xy^2 & -2y^3 \end{pmatrix}.
\end{align*}
One verifies that $\widetilde\PPP^3(\Kref)\subset \ker(\dDiv)$.
Note that this space is only defined on the reference element.
We transform it to a physical element by an appropriate transformation, discussed in the next section.

\subsection{Piola--Kirchhoff transformation}
Given a triangle $K$ with vertices $\ba_K$, $\bb_K$, $\bc_K$ or a parallelogram with vertices $\ba_K$, $\bb_K$, $\bd_K$, $\bc_K$ (ordered in positive direction), denote by $F_K\colon \Kref \to K$ the affine mapping
\begin{align*}
  \widehat\bx\mapsto \bB_K \widehat\bx + \ba_K := \begin{pmatrix} \bb_K-\ba_K & \bc_K-\ba_K\end{pmatrix}\widehat\bx + \ba_K.
\end{align*}
Note that $\det(\bB_K)>0$.
The Piola transformation $\piola_K\colon \Hdivset\Kref \to \Hdivset{K}$, $\widehat\bphi\mapsto \piola_K(\widehat\bphi)=:\bphi$ is given by
\begin{align*}
  \det(\bB_K) \bphi\circ F_K = \bB_K\widehat\bphi.
\end{align*}
We recall some of its properties.
\begin{lemma}\label{lem:piola}
  The Piola transformation is an isomorphism. Moreover, 
  \begin{align*}
   \piola_K(\widehat\bphi)\in \RT^p(K) \quad\Longleftrightarrow\quad
   \widehat\bphi\in \RT^p(\Kref).
  \end{align*}
\end{lemma}

The Piola--Kirchhoff transformation $\pk_K\colon \HdivDivset\Kref\to \HdivDivset{K}$, $\widehat\bM\mapsto \pk_K(\widehat\bM)=:\bM$ is given by
\begin{align*}
  \det(\bB_K) \bM\circ F_K = \bB_K \widehat\bM \bB_K^\top. 
\end{align*}
In the next lemma we collect some of its properties, see, e.g.,~\cite[Section~3.1]{PechsteinSchoeberl11} and~\cite[Section~4.1]{KLove2}. 
We note that~\cite{KLove2} only deals with triangular meshes. However, the proof of the next result holds verbatim for parallelogram meshes.
\begin{lemma}\label{lem:piolak}
  The Piola--Kirchhoff transformation is an isomorphism. Furthermore, 
  \begin{align*}
    \dual{\trdDiv{K}\bM}z_{\partial K} &= \dual{\trdDiv{\Kref}\widehat\bM}{\widehat z}_{\partial\Kref}, \\
    \dDiv\bM\circ F_K &= \widehat{\dDiv}\widehat\bM, \\
    \norm{\bM}K &\eqsim h_K \norm{\widehat\bM}{\Kref},
\end{align*}
for all $\bM\in\HdivDivset{K}$, $z\in H^2(K)$ with $\widehat z = z\circ F_K$, $\bM = \pk_K(\widehat\bM)$.
\end{lemma}

\section{$\HdivDivset\Omega$ elements}\label{sec:femSpace}
We describe the local finite element space for triangles in Section~\ref{sec:localSpace} and for parallelograms in Section~\ref{sec:localSpace:rect}. The global finite element space together with a canonical interpolation operator is discussed in Section~\ref{sec:globalSpace}.

\subsection{Local finite element for triangles}\label{sec:localSpace}
For triangles $K\in\TT$ we define the local spaces
\begin{align}
  \XXX_\triangle(K) = \sym(\RT^0(K)\otimes\RT^1(K)).
\end{align}

\begin{proposition}\label{propSpace}
  The following properties hold for any triangle $K\in\TT$:
\begin{enumerate}[label=(\alph*)]
  \item\label{propSpace:c} $\widehat\bM\in\XXX_\triangle(\Kref) \,\Longleftrightarrow\, \pk_K(\widehat\bM)\in \XXX_\triangle(K)$, 
  \item\label{propSpace:a} decomposition
  \begin{align*}
    \XXX_\triangle(K) = \PPPsym^1(K) \oplus \widetilde\PPP^2(K) \oplus \bx\bx^\top P^1(K),
  \end{align*}
  \item\label{propSpace:b} $\dim(\XXX_\triangle(K))=15$, 
  \item\label{propSpace:d} $\PPPsym^1(K)\subseteq \XXX_\triangle(K) \subseteq \PPPsym^3(K)$,
  \item\label{propSpace:e} $\dDiv(\XXX_\triangle(K)) = P^1(K)$,
  \item\label{propSpace:f} $\trdDiv{K,E,\normal}(\bM) \in P^1(E)$ and 
  \item\label{propSpace:g} $\trdDiv{K,E,\tangential}(\bM) \in P^1(E)$ for all $E\in\edges_K$, $\bM\in\XXX_\triangle(K)$.
\end{enumerate}
\end{proposition}
\begin{proof}
  \noindent\boxed{\textbf{\ref{propSpace:c}.}}
  This follows from the properties of the Piola transformation (Lemma~\ref{lem:piola}) and the definition of the Piola--Kirchhoff transformation.
  Let $\widehat\bM\in \XXX_\triangle(\Kref)$. Then, there exist $\widehat\bphi_j\in \RT^0(\Kref)$ and $\widehat\bpsi_j\in \RT^1(\Kref)$, $j=1,\ldots,n<\infty$, such that
  \begin{align*}
  \widehat\bM = \sum_{j=1}^{n}\frac12(\widehat\bphi_j\widehat\bpsi_j^\top+\widehat\bpsi_j\widehat\bphi_j^\top).
  \end{align*}
  Using $\bM:=\pk_K(\widehat\bM)$ we have that
  \begin{align*}
    \det(\bB_K)\bM\circ F_K = \sum_{j=1}^{n}\frac12(\bB_K\widehat\bphi_j)(\bB_K\widehat\bpsi_j)^\top + \sum_{j=1}^{n}\frac12(\bB_K\widehat\bpsi_j)(\bB_K\widehat\bphi_j)^\top.
  \end{align*}
  Using $\bphi_j:=\piola_K(\widehat\bphi_j) = \frac1{\det\bB_K} \bB_K\widehat\bphi_j \circ F_K^{-1}\in \RT^0(K)$ and $\bpsi_j:=\piola_K(\widehat\bpsi_j)\in\RT^1(K)$ in the latter identity we find that
  \begin{align*}
    \bM = \det(\bB_K) \sum_{j=1}^n \frac12 \big(\bphi_j\bpsi_j^\top + \bpsi_j\bphi_j^\top\big) \in\XXX_\triangle(K).
  \end{align*}
  This concludes the proof by noting that the Piola and Piola--Kirchhoff transformations are isomorphisms.

  \noindent\boxed{\textbf{\ref{propSpace:a}.}} This follows by simple calculations.
    For the sake of completeness we give details. Let $\widetilde{\XXX}_\triangle(K) = \PPPsym^1(K) \oplus \widetilde\PPP^2(K) \oplus \bx\bx^\top P^1(K)$. Any element in $\XXX_\triangle(K)$ can be written as a linear combination of tensors of the form
    \begin{align*}
      \frac{\bphi\bpsi^\top+\bpsi\bphi^\top}{2} + \frac{\bx\bpsi^\top+\bpsi\bx^\top}2\phi 
      +\frac{\bphi\bx^\top + \bx\bphi^\top}2\psi + \bx\bx^\top \phi\psi =: \sum_{j=1}^4 \bM_j
    \end{align*}
    with $\phi\in P^0(K)$, $\bphi\in \bP^0(K)$, $\psi\in P^1(K)$, $\bpsi\in \bP^1(K)$.
    First, we show that $\widetilde{\XXX}_\triangle(K)\subseteq \XXX_\triangle(K)$. From the last displayed formula we see that we only need to prove that $\widetilde\PPP^2(K)\subseteq \XXX_\triangle(K)$. The three elements~\eqref{eq:tildeP2span} are represented by taking $\phi=1$, $\bphi = 0$, $\psi = 0$ and $\bpsi = \bpsi_j$ ($j=1,2,3$) with $\bpsi_1 = (x,-y)^\top$, $\bpsi_2 = (y,0)^\top$, $\bpsi_3 = (0,x)^\top$.
    To prove the other inclusion $\XXX_\triangle(K)\subseteq \widetilde\XXX_\triangle(K)$ note that $\bM_1\in \PPPsym^1(K)$, $\bM_4\in \bx\bx^\top P^1(K)$ so that it only remains to prove that $\bM_2,\bM_3\in \widetilde\XXX_\triangle(K)$. Recall that $\bM_2 = \phi(\bx\bpsi^\top+\bpsi\bx^\top)/2$ with $\phi\in P^0(K)$, $\bpsi\in \bP^1(K)$. By linearity it is sufficient to prove that $\bM_2\in\widetilde\XXX_\triangle(K)$ for a basis of $\bP^1(K)$. If $\bpsi\in \bP^0(K)$ then $\bM_2 \in \PPPsym^1(K)$. 
    A basis of $P_\mathrm{hom}^1(K)^2$ is given by the vectors $\bpsi_j$ ($j=1,2,3$) and $\bpsi_4 = (x,y)^\top$.
Clearly, $\bM_2\in \bx\bx^\top P^0(K)\subset \widetilde\XXX_\triangle(K)$ for $\bpsi=\bpsi_4$. Further, $\bM_2 \in \widetilde\PPP^2(K)$ for $j=1,2,3$ which follows with similar arguments as before. 
Lastly, recall that $\bM_3 = \psi(\bphi\bx^\top + \bx\bphi^\top)/2$ with $\psi\in P^1(K)$, $\bphi\in \bP^0(K)$. If $\psi\in P^0(K)$ then $\bM_3\in \PPPsym^1(K)\subset \widetilde\XXX_\triangle(K)$. It is sufficient to check if $\bM_3\in \widetilde\XXX_\triangle(K)$ for all combinations of $\bphi = (1,0)^\top$, $\bphi =(0,1)^\top$, $\psi = x$, $\psi=y$. This gives the following four tensors
\begin{align*}
  \bM_3^{(1)} = \begin{pmatrix} x^2 &\tfrac{xy}2 \\ \tfrac{xy}2 & 0\end{pmatrix}, \,
  \bM_3^{(2)} = \begin{pmatrix} xy &\tfrac{y^2}2 \\ \tfrac{y^2}2 & 0\end{pmatrix}, \,
  \bM_3^{(3)} = \begin{pmatrix} xy &\tfrac{x^2}2 \\ \tfrac{x^2}2 & 0\end{pmatrix}, \,
  \bM_3^{(4)} = \begin{pmatrix} 0 &\tfrac{xy}2 \\ \tfrac{xy}2 & y^2\end{pmatrix}.
\end{align*}
Note that $\bM_3^{(2)},\bM_3^{(3)} \in \widetilde\PPP^2(K)$ by~\eqref{eq:tildeP2span}. The other two elements can be written as linear combinations of elements from $\widetilde\PPP^2(K)$ and $\bx\bx^\top P^1(K)$, e.g., 
\begin{align*}
  \bM_3^{(1)} = \begin{pmatrix} x^2 &\tfrac{xy}2 \\ \tfrac{xy}2 & 0\end{pmatrix}
  = \frac12 \bx\bx^\top + \frac12\begin{pmatrix} x^2 & 0\\ 0 & -y^2\end{pmatrix} \in \bx\bx^\top P^1(K)+\widetilde\PPP^2(K).
\end{align*}
Therefore, $\bM_3\in \widetilde\XXX_\triangle(K)$, thus concluding the proof of the inclusion $\XXX_\triangle(K)\subseteq \widetilde \XXX_\triangle(K)$.

  \noindent\boxed{\textbf{\ref{propSpace:b}.}}
  We use the decomposition from~\ref{propSpace:a} and count the dimensions, $\dim(\PPPsym^1(K)) = 9$, $\dim(\widetilde\PPP^2(K)) = 3$, $\dim(\bx\bx^\top P^1(K)) = 3$.

  \noindent\boxed{\textbf{\ref{propSpace:d}.}}
  Follows by definition as can be seen from decomposition~\ref{propSpace:a}.

  \noindent\boxed{\textbf{\ref{propSpace:e}.}}
  Note that $\bx\bx^\top P^1(K)\subset \XXX_\triangle(K)$ by decomposition~\ref{propSpace:a}.
  One easily verifies that $\dDiv\colon \bx\bx^\top P^1(K)\to P^1(K)$ is surjective
  by calculating $\dDiv (\bx\bx^\top \phi)=6\phi+3\bx^\top\grad\phi$ for
  $\phi\in P^1(K)$, and comparing coefficients.

  \noindent\boxed{\textbf{\ref{propSpace:f}.}}
  Let $\bphi\in \RT^0(K)$ and $\bpsi\in\RT^1(K)$
  and consider
  \begin{align*}
    \bM = \frac12(\bphi\bpsi^\top + \bpsi\bphi^\top) \in\XXX_\triangle(K).
  \end{align*}
  Recall that $\bphi\cdot\normal|_E\in P^0(E)$, $\bpsi\cdot\normal|_E\in P^1(E)$ for any $E\in\edges_K$. Therefore,
  \begin{align*}
    \trdDiv{K,E,\normal}(\bM) = \normal\cdot\bM\normal|_E = \frac12\big( (\bphi\cdot\normal)(\bpsi\cdot\normal) + (\bpsi\cdot\normal)(\bphi\cdot\normal) \big)|_E \in P^1(E).
  \end{align*}
  The assertion follows for all $\bM\in\XXX_\triangle(K)$ by linearity of the trace operator.

  \noindent\boxed{\textbf{\ref{propSpace:g}.}}
  Let $\bM\in\XXX_\triangle(K)$.
  Recall that
  \begin{align*}
    \trdDiv{K,E,\tangential}(\bM) = (\normal\cdot\Div\bM+\partial_{\tangential}(\tangential\cdot\bM\normal))|_E.
  \end{align*}
  We write 
  \begin{align*}
    \bM = \bM_1 + \bM_2
  \end{align*}
  with $\bM_1\in \PPPsym^1(K)+\widetilde\PPP^2(K)$ and $\bM_2=\bx\bx^\top\eta$ for some $\eta\in P^1(K)$.
  Looking at the polynomial degrees it is clear that $\trdDiv{K,E,\tangential}(\bM_1) \in P^1(E)$.
  A simple computation shows 
  \begin{align*}
    \Div\bM_2 = \bx(\underbrace{\bx^\top\nabla\eta+3\eta}_{\in P^1(K)}).
  \end{align*}
  Recall that $\bx\cdot\normal|_E \in P^0(E)$. Therefore, 
  \begin{align*}
    \normal\cdot\Div\bM_2|_E \in P^1(E)
    \quad\text{and}\quad 
    \tangential\cdot\bM_2\normal|_E \in P^2(E).
  \end{align*}
  Consequently, $\partial_\tangential(\tangential\cdot\bM_2\normal)|_E\in P^1(E)$ which proves the assertion.
\end{proof}

We define the following degrees of freedom for space $\XXX_\triangle(K)$,
\begin{subequations}\label{eq:dofs}
  \begin{alignat}{3}
    &\text{the moments }&&\frac{1}{\norm{\ell_{E,k}}E^2}\dual{\trdDiv{K,E,\normal}(\bM)}{\ell_{E,k}}_E &\quad k=0,1, \forall E\in\edges_K, \label{eq:dofs:nMn}\\
  &\text{the moments }&&\dual{\trdDiv{K,E,\tangential}(\bM)}{\ell_{E,k}}_E &\quad k=0,1, \forall E\in\edges_K,\label{eq:dofs:tMn} \\
  &\text{the values } &&\jump{\tangential\cdot\bM\normal}_{\partial K}(z) &\quad\forall z\in \vertices_K. \label{eq:dofs:jumps}
\end{alignat}
\end{subequations}
The scaling factor $1/\norm{\ell_{E,k}}E^2$ in~\eqref{eq:dofs:nMn} is used to ensure that all elements of the basis which are dual to the degrees of freedom have the same scaling, see Appendix~\ref{sec:appendix} below.

\begin{theorem}\label{thm:dofs}
  Degrees of freedom~\eqref{eq:dofs} are unisolvent for space $\XXX_\triangle(K)$.
\end{theorem}
\begin{proof}
  We have 15 degrees of freedom and $\dim(\XXX_\triangle(K))=15$. Thus, it suffices to prove that if for $\bM\in\XXX_\triangle(K)$ all degrees of freedom~\eqref{eq:dofs} vanish, then $\bM=0$.
  First, observe that by Proposition~\ref{propSpace}, the vanishing degrees of freedom, and the definition of trace operator $\trdDiv{K}$, it follows from \eqref{eq:jump} that
  \begin{align*}
    0 = \trdDiv{K}\bM,
  \end{align*}
  By Lemma~\ref{lem:piolak} it also follows that $\trdDiv{\Kref}{\widehat\bM} =0$ where $\widehat\bM = \pk_K^{-1}\bM \in \XXX_\triangle(\Kref)$.
  W.l.o.g. we can thus assume that $K=\Kref$ for the remainder of the proof.

  Note that $u:=\dDiv\bM\in P^1(K)$. Then, integration by parts proves
  \begin{align*}
    \norm{\dDiv\bM}{K}^2 = \ip{\dDiv\bM}u_K = \ip{\bM}{\Ggrad u}_K + \dual{\trdDiv{K}(\bM)}{u|_{\partial K}}_{\partial K} = 0. 
  \end{align*}
  By the decomposition of Proposition~\ref{propSpace} this means that $\bM=\bM_1+\bM_2\in \PPPsym^1(K)\oplus\widetilde\PPP^2(K)$.
  Note that $\trdDiv{K,E,\tangential}(\bM_1)\in P^0(E)$, hence, $\dual{\trdDiv{K,E,\tangential}(\bM_1)}{\ell_{E,1}}_E = 0$ for all $E\in\edges_K$ (remember that $\ell_{E,1}$ is the Legendre
polynomial of degree 1, transformed to $E$).
  From 
  \begin{align*}
    \dual{\trdDiv{K,E,\tangential}(\bM_2)}{\ell_{E,1}}_E = \dual{\trdDiv{K,E,\tangential}(\bM)}{\ell_{E,1}}_E = 0 \quad\forall E\in\edges_K
  \end{align*}
  one finds after a short calculation that $\bM_2=0$, thus, $\bM = \bM_1\in \PPPsym^1(K)$. 

  It remains to show that $\bM_1=0$. To that end, note that a basis of $\PPPsym^1(K)$ is given by
  \begin{align*}
    \bN_j = \sym(\tangential_{j+1}\tangential_{j+2}^\top), \quad \bN_{j+3} = (\eta_{j+1}-\eta_j)\bN_j, 
    \quad \bN_{j+6} = \eta_{j+2}\bN_j,\quad  j=1,2,3,
  \end{align*}
  Here, $\eta_j$, $j=1,2,3$, are the barycentric coordinates of $K$, numbered in
cyclic form, and $E_j$, $j=1,2,3$, denote the three edges of $K$ connecting the nodes $z_j$ and $z_{j+1}$. Their tangential vectors are $\tangential_j$ and we use cyclic indexing, i.e., $z_{3+j} = z_j$ for $j=1,2,3$. 
  Observe that $\normal\cdot\bN_{j}\normal|_{\partial K} = 0$ for $j=7,8,9$ and $\normal\cdot\bN_{j}\normal|_{E_k} =c_j\delta_{jk}$ for $j,k=1,2,3$  and $\normal\cdot\bN_{3+j}\normal|_{E_k} =d_j\delta_{jk}\ell_{E_j,1}$ (with non-vanishing constants $c_j$, $d_j$).
  From $\dual{\trdDiv{K,E,\normal}(\bM)}{p}_E=0$ for all $p\in P^1(E)$ and $E\in\edges_K$ we conclude
  that $\bM = \sum_{j=7}^9 \alpha_j \bN_j$.
  Then, a simple computation shows that $\jump{\tangential\cdot\bN_{6+j}\normal}|_{\partial K}(z_{k}) = e_j\delta_{j+2,k}$ with non-vanishing constant $e_j$.
  Thus, $\jump{\tangential\cdot\bM\normal}|_{\partial K}(z)=0$ for all $z\in \vertices_K$ shows $\bM = 0$ and finishes the proof.
\end{proof}

\begin{remark} \label{rem:dofs}
In Section~\ref{sec:appendix:triangle}
we construct a basis of $\XXX_\triangle(\Kref)$ that is dual to degrees
of freedom \eqref{eq:dofs}. Of course, this also proves their unisolvency.
\end{remark}

\subsection{Local finite element for parallelograms}\label{sec:localSpace:rect}
Throughout this section, $K = F_K(\Kref)$ with $\Kref=\Kref_\square$.
The lowest-order Raviart--Thomas space on the reference square is
\begin{align*}
  \bQ^0(\Kref) = \bP^0(\Kref) \oplus \linhull\left\lbrace (x,y)\mapsto \begin{pmatrix}x\\0\end{pmatrix},
  \, (x,y)\mapsto \begin{pmatrix}0\\y\end{pmatrix}
  \right\rbrace.
\end{align*}
On the physical element we use the Piola transformation and define
\begin{align*}
  \bQ^0(K) = \piola_K(\bQ^0(\Kref)).
\end{align*}
The idea is to define the local space in the same spirit as before for triangles. Set
\begin{align*}
  \XXX_\square(\Kref) = \sym(\bQ^0(\Kref)\otimes \RT^1(\Kref))
 \quad\text{and}\quad
  \XXX_\square(K) =  \sym(\bQ^0(K)\otimes \RT^1(K)).
\end{align*}
Noting that $\RT^0(K)\subset \bQ^0(K)$ we see that $\XXX_\triangle(K)\subset \XXX_\square(K)$. Here, we stress that $\XXX_\triangle(K)$ is well defined for $K$ being a parallelogram.

\begin{proposition}\label{propSpace:rect}
  The following properties hold for any parallelogram $K\in\TT$:
\begin{enumerate}[label=(\alph*)]
  \item\label{propSpace:rect:c} $\widehat\bM\in\XXX_\square(\Kref) \,\Longleftrightarrow\, \pk_K(\widehat\bM)\in \XXX_\square(K)$, 
  \item\label{propSpace:rect:a} decomposition
  \begin{align*}
    \XXX_\square(K) = \XXX_\triangle(K) \oplus \pk_K(\widetilde\PPP^3(\Kref)) = \PPPsym^1(K) \oplus \widetilde\PPP^2(K) \oplus \pk_K(\widetilde\PPP^3(\Kref)) \oplus \bx\bx^\top P^1(K),
  \end{align*}
  \item\label{propSpace:rect:b} $\dim(\XXX_\triangle(K))=20$, 
  \item\label{propSpace:rect:d} $\PPPsym^1(K)\subseteq \XXX_\triangle(K) \subseteq \XXX_\square(K) \subseteq \PPPsym^3(K)$,
  \item\label{propSpace:rect:e} $\dDiv(\XXX_\square(K)) = P^1(K)$,
  \item\label{propSpace:rect:f} $\trdDiv{K,E,\normal}(\bM) \in P^1(E)$ and 
  \item\label{propSpace:rect:g} $\trdDiv{K,E,\tangential}(\bM) \in P^1(E)$ for all $E\in\edges_K$, $\bM\in\XXX_\triangle(K)$.
\end{enumerate}
\end{proposition}
\begin{proof}
  \noindent\boxed{\textbf{\ref{propSpace:rect:c}.}} Proof is identical to the corresponding one from Proposition~\ref{propSpace}. We therefore skip the details.

  \noindent\boxed{\textbf{\ref{propSpace:rect:a}.}} 
  This can be seen by noting that a straightforward computation on the reference element proves $\XXX_\square(\Kref) = \XXX_\triangle(\Kref) \oplus \widetilde\PPP^3(\Kref)$ and that $\XXX_\triangle(K) = \pk_K(\XXX_\triangle(\Kref))$.
  
  \noindent\boxed{\textbf{\ref{propSpace:rect:b}.}} Follows by counting dimensions of $\XXX_\square(\Kref)$ or in decomposition~\ref{propSpace:rect:a}.

  \noindent\boxed{\textbf{\ref{propSpace:rect:d}.}} 
  This follows directly by definition and Proposition~\ref{propSpace}.

  \noindent\boxed{\textbf{\ref{propSpace:rect:e}.}} 
  This follows from $\dDiv(\XXX_\triangle(K)) = P^1(K)$, $\XXX_\triangle(K)\subset \XXX_\square(K)$ and $\dDiv(\XXX_\square(K))\subseteq P^1(K)$.

  \noindent\boxed{\textbf{\ref{propSpace:rect:f}.}} 
  We argue as in the proof of Proposition~\ref{propSpace}. Let $\bphi\in \bQ^0(K)$, $\bpsi\in \RT^1(K)$. Note that $\bphi\cdot\normal|_E\in P^0(E)$ and $\bpsi\cdot\normal|_E\in P^1(E)$ for any $E\in\edges_K$.
  It follows that for $\bM = \bphi\bpsi^\top + \bpsi\bphi^\top\in \XXX_\square(K)$ we have
  \begin{align*}
    \normal\cdot\bM\normal|_E = (\bphi\cdot\normal)(\bpsi\cdot\normal)|_E + (\bpsi\cdot\normal)(\bphi\cdot\normal)|_E \in P^1(E).
  \end{align*}
  Recall that any $\bM\in \XXX_\square(K)$ can be written as $\bM = \sum_{j=1}^n \sym(\bphi_j\bpsi_j^\top)$ for some 
  $\bphi_j\in\bQ^0(K)$, $\bpsi_j\in\RT^1(K)$, and $n\in\N$.
  With the previous observation and linearity we conclude that
  \begin{align*}
    \normal\cdot\bM\normal|_E \in P^1(E).
  \end{align*}

  \noindent\boxed{\textbf{\ref{propSpace:rect:g}.}} 
  Let $E\in \edges_K$ be fixed and $\bphi\in \bQ^0(K)$, $\bpsi\in\RT^1(K)$ be given. Set $\bM = \sym(\bphi\bpsi^\top)$ and observe that
  \begin{align*}
    \tangential\cdot \bM\normal|_E\in P^2(E) 
  \end{align*}
  and, therefore, $\partial_\tangential(\tangential\cdot\bM\normal)|_E \in P^1(E)$.
  By linearity we conclude that $\partial_\tangential(\tangential\cdot\bM\normal)|_E \in P^1(E)$ for all $\bM\in\XXX_\square(K)$.

  Some tensor calculus yields
  \begin{align*}
    \det(\bB_K)\Div\bM \circ F_K = \bB_K \widehat\Div\widehat\bM.
  \end{align*}
Multiplying with the normal on $E$, denoted by $\normal$, one finds that $(B_K^{-1})^\top\normal=c\widehat\normal$ and
  \begin{align*}
    \det(\bB_K)\normal\cdot\Div\bM \circ F_K = c \widehat\normal\cdot\widehat\Div\widehat\bM
  \end{align*}
  with $c$ being a non-vanishing constant. The proof is finished if $\widehat\normal\cdot\widehat\Div\widehat\bM\in P^1(\Eref)$ on edge $\Eref = F_K^{-1}(E)$.
  If $\widehat\bM\in \XXX_\triangle(\Kref)$ then we have already seen this in the proof for triangles (there we have shown this directly). The same argumentation applies here.
  W.l.o.g. we can thus assume that $\widehat\bM\in \widetilde\PPP^3(\Kref)$. There are only two basis function with degree $3$ polynomials. For these two we find
  \begin{align*}
    \Div\begin{pmatrix} 2x^3 & -x^2y \\ -x^2y & -4xy^2 \end{pmatrix} = \begin{pmatrix} 5x^2\\-10xy\end{pmatrix},
    \quad
    \Div\begin{pmatrix} 4x^2y & xy^2 \\ xy^2 & -2y^3 \end{pmatrix} = \begin{pmatrix} 10xy \\ -5y^2\end{pmatrix}.
  \end{align*}
  One verifies that the normal traces of these two vectors are polynomials of degree $\leq 1$ on each edge $\Eref\in\edges_{\Kref}$. This concludes the proof.
\end{proof}

We define the same degrees of freedom~\eqref{eq:dofs} for space $\XXX_\square(K)$ with $K$ being a parallelogram as for the triangular element, of course noting that we now have four vertices and edges
instead of three.
% \begin{subequations}\label{eq:dofs:rect}
% \begin{align}
%   \text{the moments }\dual{\trdDiv{K,E,\normal}(\bM)}{\ell_{E,k}}_E &\quad k=0,1, \forall E\in\edges_K, \label{eq:dofs:rect:nMn}\\
%   \text{the moments }\dual{\trdDiv{K,E,\tangential}(\bM)}{\ell_{E,k}}_E &\quad k=0,1, \forall E\in\edges_K,\label{eq:dofs:rect:tMn} \\
%   \text{the values }\jump{\tangential\cdot\bM\normal}_{\partial K}(z) &\quad\forall z\in \vertices_K. \label{eq:dofs:rect:jumps}
% \end{align}
% \end{subequations}

\begin{theorem}\label{thm:dofs:rect}
  Degrees of freedom~\eqref{eq:dofs} are unisolvent for space $\XXX_\square(K)$ with $K$ being a parallelogram.
\end{theorem}
\begin{proof}
  We have 20 degrees of freedom and $\dim(\XXX_\square(K))=20$. Thus, it suffices to prove that if for $\bM\in\XXX_\square(K)$ all degrees of freedom~\eqref{eq:dofs} vanish, then $\bM=0$.
  For the remainder of the proof we argue similarly as in the proof of Theorem~\ref{thm:dofs}.
  First, observe that by Proposition~\ref{propSpace}, the vanishing degrees of freedom, and the definition of trace operator $\trdDiv{K}$, it follows that
  \begin{align*}
    0 = \trdDiv{K}\bM.
  \end{align*}
  By Lemma~\ref{lem:piolak} it also follows that $\trdDiv{\Kref}{\widehat\bM} =0$ where $\widehat\bM = \pk_K^{-1}\bM \in \XXX_\triangle(\Kref)$.
  W.l.o.g. we can thus assume that $K=\Kref=\Kref_\square$ for the remainder of the proof.

  With $u:=\dDiv\bM\in P^1(K)$ integration by parts proves
  \begin{align*}
    \norm{\dDiv\bM}{K}^2 = \ip{\dDiv\bM}u_K = \ip{\bM}{\Ggrad u}_K + \dual{\trdDiv{K}(\bM)}{u|_{\partial K}}_{\partial K} = 0. 
  \end{align*}
  By Proposition~\ref{propSpace:rect} this means that $\bM\in\PPPsym^1(K)\oplus\widetilde\PPP^2(K)\oplus\widetilde\PPP^3(K)$.
  
  The condition $\trdDiv{K,E,\tangential}\bM=0$ for all $E\in\edges_K$ implies that 
  \begin{align*}
    \bM \in \PPPsym^0(K) \oplus \symCurl\big(\linhull\set{\bphi_j}{j=1,\ldots,6}\big)
  \end{align*}
  where
  \begin{alignat*}{3}
    \bphi_1(x,y) &= \begin{pmatrix} 0 \\ x^2\end{pmatrix}, \quad \bphi_2(x,y) &\,=\,& \begin{pmatrix} y^2 \\ 0\end{pmatrix}, \quad &\bphi_3(x,y) =&\,  \begin{pmatrix} xy \\ 0\end{pmatrix}, \\
    \bphi_4(x,y) &= \begin{pmatrix} 0\\xy\end{pmatrix}, \quad \bphi_5(x,y) &\,=\,& \begin{pmatrix} xy^2 \\ 0\end{pmatrix}, \quad &\bphi_6(x,y) =&\, \begin{pmatrix} 0 \\ x^2 y\end{pmatrix}.
  \end{alignat*}
  The latter condition requires some tedious but simple computations which are not shown here for simplicity.
  The condition $\trdDiv{K,E,\normal}\bM = 0$ then yields
  \begin{align*}
    \bM \in \linhull\left\lbrace (x,y)\mapsto \begin{pmatrix} 0 & 1\\ 1 & 0\end{pmatrix} \right\rbrace. 
  \end{align*}
  Finally, from $\jump{\tangential\cdot\bM\normal}_{\partial K}(z) = 0$ for $z=(0,0)$ one concludes
  that $\bM = 0$, finishing the proof.
\end{proof}

\begin{remark}
Remark~\ref{rem:dofs} applies in this case as well. A dual basis is given in
Section~\ref{sec:appendix:square}.
\end{remark}

\subsection{Global space and canonical interpolation}\label{sec:globalSpace}
We generically set
\begin{align*}
  \XXX(K) = \begin{cases} 
    \XXX_\triangle(K) & \text{if } K \text{ is a triangle}, \\
    \XXX_\square(K) & \text{if } K \text{ is a parallelogram}.
  \end{cases}
\end{align*}
The global space is given by
\begin{align*}
  \XXX(\TT) := \set{\bM\in\LLLts(\Omega)}{\bM|_K\in \XXX(K), \, K\in\TT}\cap\HdivDivset\Omega.
\end{align*}
The $\HdivDivset\Omega$-conformity requires to assign unique edge-degrees of freedom
\eqref{eq:dofs:nMn}, \eqref{eq:dofs:tMn} ($4$ per edge),
and vertex-degrees of freedom \eqref{eq:dofs:jumps}
($3$ per triangles and $4$ per parallelogram),
subject to the jump constraints at interior nodes of the mesh,
cf.~Proposition~\ref{prop_conf}. Therefore, the dimension of this space is
\begin{align*}
  \dim(\XXX(\TT)) = 4\#\edges + 3\#\TT_\triangle + 4\#\TT_\square - \#\vertices_\Omega.
\end{align*}
Here, $\TT_\triangle\subseteq \TT$ and $\TT_\square\subseteq\TT$ denote the sets of all triangles and parallelograms in $\TT$, respectively.

We define the canonical interpolation operator 
\begin{align*}
  \projHdDiv_\TT\colon \HdivDivset\Omega\cap \HHH^{1+r}(\TT)\to \XXX(\TT) \quad (r>1/2)
\end{align*}
by interpolation in the degrees of freedom, i.e., 
\begin{align}
\begin{split}\label{eq:defprojHdDiv}
  \frac1{\norm{\ell_{E,k}}E^2}\dual{\trdDiv{K,E,\normal}(\bM-\projHdDiv_\TT\bM)}{\ell_{E,k}}_E &= 0,\\
  \dual{\trdDiv{K,E,\tangential}(\bM-\projHdDiv_\TT\bM)}{\ell_{E,k}}_E &= 0,\\
  \jump{\tangential\cdot(\bM-\projHdDiv_\TT\bM)\normal}_{\partial K}(z) &= 0
\end{split}
\end{align}
for $k=0,1$ and all $E\in\edges_K$, $z\in \vertices_K$ and $K\in\TT$.

\begin{proposition}\label{prop:proj}
  Operator $\projHdDiv_\TT$ is well defined for $r>1/2$ and is a projection. It has
  the commutativity property
  \begin{align*}
    \dDiv\circ\projHdDiv_\TT = \Pi_\TT^1\circ\dDiv
  \end{align*}
  and approximation property
  \begin{align*}
    \norm{\bM-\projHdDiv_\TT\bM}\Omega \lesssim \norm{h_\TT}{L_\infty(\Omega)}^{\min\{1+r,2\}}\norm{\bM}{H^{\min\{1+r,2\}}(\Omega)}
  \end{align*}
  for all $\bM\in \HHH^{\min\{1+r,2\}}(\Omega)\cap\HdivDivset\Omega$ with $r>1/2$.
\end{proposition}
\begin{proof}
  Let $\bM\in \HHH^{1+r}(\Omega)\cap\HdivDivset\Omega$ with $r>1/2$. Then it is clear that $\trdDiv{K,E,\normal}\bM\in L^2(E)$, $\trdDiv{K,E,\tangential}\bM\in L^2(E)$, $\jump{\tangential\cdot\bM\normal}_{\partial K}(z)\in\R$.
  Therefore, operator $\projHdDiv_\TT\bM|_K$ is well defined for each $K\in\TT$. By definition we also conclude that $\projHdDiv_\TT\bM\in \XXX(\TT)$ since for an interior node $z\in\vertices_\Omega$ we observe
  \begin{align*}
    \sum_{K\in\patch_z} \jump{\tangential\cdot\projHdDiv_\TT\bM\normal}_{\partial K}(z) = \sum_{K\in\patch_z} \jump{\tangential\cdot\bM\normal}_{\partial K}(z) = 0.
  \end{align*}
  To see the commutativity property, let $K\in\TT$ and $u\in P^1(K)$ be given. Integration by parts and the definition of $\projHdDiv_\TT$, see~\eqref{eq:defprojHdDiv}, show that
  \begin{align*}
    \ip{\dDiv\projHdDiv_\TT\bM}{u}_K &= \ip{\projHdDiv_\TT\bM}{\Ggrad u}_K + \dual{\trdDiv{K}\projHdDiv_\TT\bM}{u}_{\partial K} 
    %\\
    %&= \sum_{E\in\edges_K} \dual{\trdDiv{K}\projHdDiv_\TT\bM}{u}_E 
    \\ 
    &= \sum_{E\in\edges_K} \big(\dual{\trdDiv{K,E,\tangential}\projHdDiv_\TT\bM}{u}_E - \dual{\trdDiv{K,E,\normal}\projHdDiv_\TT\bM}{\partial_\normal u}_E\big) 
    \\
    &\qquad - \sum_{z\in\vertices_K} \jump{\tangential\cdot\projHdDiv_\TT\bM\normal}_{\partial K}(z)u(z) 
  \\
  &= \sum_{E\in\edges_K} \big(\dual{\trdDiv{K,E,\tangential}\bM}{u}_E - \dual{\trdDiv{K,E,\normal}\bM}{\partial_\normal u}_E\big) 
    \\
    &\qquad - \sum_{z\in\vertices_K} \jump{\tangential\cdot\bM\normal}_{\partial K}(z)u(z)
  %\\
  %  &= \sum_{E\in\edges_K} \dual{\trdDiv{K}\bM}{u}_E 
  \\
  &= \ip{\bM}{\Ggrad u}_K + \dual{\trdDiv{K}\bM}{u}_{\partial K} = \ip{\dDiv\bM}{u}_K.
  \end{align*}
  This proves the claim.

  Let $K\in\TT$ and $d_j\colon \XXX(K)\to\R$, $j=1,\dots,\dim(\XXX(K))$, denote the local degrees of freedom ordered as in~\eqref{eq:defprojHdDiv} and let $\bN_1,\dots,\bN_{\dim(\XXX(K))}$ denote the basis of $\XXX(K)$ with $d_j(\bN_k) = \delta_{jk}$, see Appendix~\ref{sec:appendix} for an explicit construction of the basis. Operator $\projHdDiv_\TT$ restricted to $K$ (and denoted by $\projHdDiv_K$) has the representation
  \begin{align*}
    \projHdDiv_K\bM = \sum_{j=1}^{\dim(\XXX(K))} d_j(\bM) \bN_j.
  \end{align*}
  If $\bM|_K\in\XXX(K)$ one concludes that $\projHdDiv_K\bM = \bM|_K$ which means that $\projHdDiv_\TT$ is a projection.

  It remains to show boundedness.
  Note that $\projHdDiv_K$ is a local projection and preserves $\PPPsym^1(K)\subset \XXX(K)$. We have that
  \begin{align*}
    \norm{(1-\projHdDiv_\TT)\bM}K = \norm{(1-\projHdDiv_\TT)(\bM-\bN)}K \leq \norm{\bM-\bN}K + \norm{\projHdDiv_\TT(\bM-\bN)}K
  \end{align*}
  for any $\bN\in\PPPsym^1(K)$. 
  For the remainder let $\bN = \Pi_K^1\bM$. Using the above representation we infer that
  \begin{align}\label{eq:proj:proof:a}
    \norm{\projHdDiv_K(\bM-\bN)}K \leq \sum_{j=1}^{\dim(\XXX(K))} |d_j(\bM-\bN)|\norm{\bN_j}K \eqsim\sum_{j=1}^{\dim(\XXX(K))} h_K|d_j(\bM-\bN)|.
  \end{align}
  The last equivalence follows from scaling properties of the basis functions, see Appendix~\ref{sec:appendix}. The proof is finished if we can show that $|d_j(\bM-\bN)| \lesssim h_K^r\norm{\bM}{H^{1+r}(K)}$ for all $j=1,\ldots,\dim(\XXX(K))$.
  Suppose that number $j$ corresponds to an edge $E$ and degree of freedom~\eqref{eq:dofs:nMn}, then
  \begin{align*}
    |d_j(\bM-\bN)| &= \frac{1}{\norm{\ell_{E,j}}E^2}|\dual{\normal\cdot(\bM-\bN)\normal}{\ell_{E,k}}_E| 
    \\
    &\lesssim h_E^{-1/2}(h_E^{1/2}\norm{\Grad(\bM-\bN)}K + h_E^{-1/2}\norm{\bM-\bN}K) 
    \lesssim h_K^{r}\norm{\bM}{H^{1+r}(K)}.
  \end{align*}
  For the latter estimates we have used the trace inequality and the approximation properties of $\Pi_K^1$.
  For indices corresponding to the other degrees of freedom one argues in a similar way.
  For instance, let $j$ refer to a degree of freedom associated with~\eqref{eq:dofs:jumps}. Using the embedding $H^{1+r}(\Kref)\to C(\overline{\Kref})$ and a scaling argument we find that
  \begin{align*}
    |d_j(\bM-\bN)|\lesssim |(\bM-\bN)(z)| &\lesssim \frac{1}{|K|^{1/2}}(\norm{\bM-\bN}K + h_K^{1+r}|\bM-\bN|_{H^{1+r}(K)})
    \\
    &\lesssim \frac{h_K^{1+r}}{|K|^{1/2}} \norm{\bM}{H^{1+r}(K)} \eqsim h_K^r \norm{\bM}{H^{1+r}(K)}.
  \end{align*}
  A similar result is found for indices corresponding to~\eqref{eq:dofs:tMn} which can be seen as follows: 
  Let $j$ correspond to one of the degree of freedom~\eqref{eq:dofs:tMn}. Then, 
  \begin{align*}
    |d_j(\bM-\bN)| &= |\dual{\normal\cdot\Div(\bM-\bN)+\partial_\tangential(\tangential\cdot(\bM-\bN)\normal)}{\ell_{E,k}}_E|
    \\
    &\lesssim \norm{\normal\cdot\Div(\bM-\bN)+\partial_\tangential(\tangential\cdot(\bM-\bN)\normal)}E \norm{\ell_{E,k}}E
    \\
    &\lesssim h_E^{r-1/2}\norm{\bM}{H^{1+r}(K)} h_E^{1/2} \eqsim  h_K^{r}\norm{\bM}{H^{1+r}(K)}.
  \end{align*}
  Putting all the estimates for $|d_j(\bM-\bN)|$ together this proves with~\eqref{eq:proj:proof:a} the estimate 
  \begin{align*}
    \norm{(1-\projHdDiv_K)\bM}K \lesssim h_K^{1+r}\norm{\bM}{H^{1+r}(K)} + \sum_{j=1}^{\dim(\XXX(K))} h_K|d_j(\bM-\bN)| \lesssim h_K^{1+r} \norm{\bM}{H^{1+r}(K)}.
\end{align*}
  Summing over all $K\in\TT$ finishes the proof.
\end{proof}

\section{Mixed finite element method}\label{sec:mixed}
As an application of our new finite element, let us consider the Kirchhoff--Love plate bending problem 
\begin{subequations}\label{eq:biharm}
\begin{align}
  \dDiv\Ctens \Ggrad u &= f \quad\text{in }\Omega, \\
  \trGgrad{\Omega}u &= \widehat\bu,
\end{align}
\end{subequations}
where $\Gamma = \partial\Omega$ and $\widehat\bu\in \trGgrad\Omega(H^2(\Omega))$.
Note that the boundary condition is often written in the form $u|_\Gamma = g_1$, $\partial_\normal u|_\Gamma  = g_2$.
Here, we assume that $\Ctens\colon \LLLts(\Omega)\to \LLLts(\Omega)$ denotes a positive definite isomorphism.
This implies that $\Ctens^{-1}$ is a positive definite isomorphism.

Introducing bending moments $\bM = \Ctens \Ggrad u$ we consider the variational mixed form:
Find $(\bM,u)\in \HdivDivset\Omega\times L_2(\Omega)$ such that
\begin{subequations}\label{eq:varform}
  \begin{alignat}{2}
    &\ip{\Ctens^{-1}\bM}{\bN}_\Omega  - \ip{u}{\dDiv\bN}_\Omega &\,=\,& -\dual{\widehat\bu}{\bN}_\Gamma, \\
    &\ip{\dDiv\bM}{v}_\Omega &\,=\,& \ip{f}v_\Omega
  \end{alignat}
  for all $(\bN,v)\in \HdivDivset\Omega\times L_2(\Omega)$.
\end{subequations}
This formulation is obtained by testing equation $\Ctens^{-1}\bM-\Ggrad u = 0$
and applying trace operator $\trGgrad{\Omega}$.

In the following, $\norm{\cdot}{3/2,1/2,\Gamma}$ is the induced trace norm of $\trGgrad{\Omega}(H^2(\Omega))$.

\begin{proposition}
  Problem~\eqref{eq:varform} is well posed. The unique solution $(\bM,u)$ satisfies $\bM=\Ctens\Ggrad u$ and $u\in H^2(\Omega)$ solves~\eqref{eq:biharm}. Furthermore, 
  \begin{align*}
    \norm{\bM}{\HdivDivset\Omega} + \norm{u}{\Omega} \lesssim \norm{f}\Omega + \norm{\widehat\bu}{3/2,1/2,\Gamma}.
  \end{align*}
\end{proposition}
\begin{proof}
  The statement follows by standard arguments from the Babu\v{s}ka--Brezzi theory, cf.~\cite{BoffiBrezziFortin,Gatica14}. 
  The right-hand side linear forms are bounded by
  definition of trace operator $\trGgrad\Omega$ and the involved norms, and the Cauchy--Schwarz inequality
  applied to $\ip{f}v_\Omega$.
  Given $\bM\in\HdivDivset\Omega$ with $\ip{\dDiv\bM}{v}_\Omega=0$ for all $v\in L_2(\Omega)$ it follows that $\dDiv\bM=0$ so that coercivity $\ip{\Ctens^{-1}\bM}{\bM}_\Omega \eqsim \norm{\bM}{\HdivDivset\Omega}^2$ holds by assumption on $\Ctens$. It only remains to note the surjectivity of 
  \begin{align*}
    \dDiv\colon \HdivDivset\Omega \to L_2(\Omega).
  \end{align*}
  In fact, given $f\in L_2(\Omega)$, we define $\bM := \Ggrad z$ where $z\in H_0^2(\Omega)$ solves $\Delta^2 z = f$, cf.~\cite{BlumRannacher80}. It follows that $\bM\in \LLLts(\Omega)$ and $\dDiv\bM = f$ as wanted. 
\end{proof}

Our mixed finite element method consists in replacing $\HdivDivset\Omega$ with $\XXX(\TT)$
and $L_2(\Omega)$ with $P^1(\TT)$, yielding: Find $(\bM_\TT,u_\TT)\in \XXX(\TT)\times P^1(\TT)$ such that
\begin{subequations}\label{eq:mixed}
  \begin{alignat}{2}
    &\ip{\Ctens^{-1}\bM_\TT}{\bN}_\Omega  - \ip{u_\TT}{\dDiv\bN}_\Omega &\,=\,& -\dual{\widehat\bu}{\bN}_\Gamma, \label{eq:mixed:a}\\
    &\ip{\dDiv\bM_\TT}{v}_\Omega &\,=\,& \ip{f}v_\Omega \label{eq:mixed:b}
  \end{alignat}
  for all $(\bN,v)\in \XXX(\TT)\times P^1(\TT)$.
\end{subequations}

The next theorem is the main result of this section.
\begin{theorem}\label{thm:mixed}
  Scheme~\eqref{eq:mixed} is well posed. Let $(\bM_\TT,u_\TT)$ denote the unique solution to~\eqref{eq:mixed} and $(\bM,u)$ the unique solution of~\eqref{eq:varform}. Then,
  \begin{align*}
    \norm{\bM-\bM_\TT}{\HdivDivset\Omega} &\lesssim \min_{\bN\in\XXX(\TT)} \norm{\bM-\bN}{\HdivDivset\Omega}, \\
    \norm{u-u_\TT}{\Omega} &\lesssim \min_{\bN\in\XXX(\TT)} \norm{\bM-\bN}{\HdivDivset\Omega} + \min_{v\in P^1(\TT)} \norm{u-v}\Omega.
  \end{align*}
  Furthermore,
  \begin{align*}
    \norm{\Ctens^{-1/2}(\bM-\bM_\TT)}{\Omega} = \min_{\bN\in\XXX(\TT), \dDiv\bN = \Pi_\TT^1\dDiv\bM} \norm{\Ctens^{-1/2}(\bM-\bN)}\Omega.
  \end{align*}
\end{theorem}
\begin{proof}
  The proof follows the usual proofs for mixed finite element methods, see, e.g.,~\cite{BoffiBrezziFortin}. 
  To see the discrete $\inf$--$\sup$ condition, let $u\in P^1(\TT)$ be given. Define $z\in H_0^1(\Omega)$ as the solution to $\Delta z = u$.
At least since Kondrat'ev \cite{Kondratiev_67_BVP} it is well known that $z\in H^{1+r}(\Omega)$ for some $1/2<r\leq 1$ depending only on $\Omega$,
see, e.g., \cite[Theorem 2]{Grisvard_76_BSE}.

 Then, by Proposition~\ref{prop:proj}, $\bN:=\projHdDiv_\TT(z\bI)$ is well defined with $\dDiv\bN = u$ and 
  \begin{align*}
    \norm{\bN}{\HdivDivset\Omega} \lesssim \norm{u}\Omega.
  \end{align*}
  It follows that
  \begin{align*}
    \sup_{\bQ\in\XXX(\TT)\setminus\{0\}} \frac{\ip{\dDiv\bQ}u_\Omega}{\norm{\bQ}{\HdivDivset\Omega}} \geq\frac{\ip{\dDiv\bN}u_\Omega}{\norm{\bN}{\HdivDivset\Omega}} 
    \gtrsim \norm{u}\Omega.
  \end{align*}
  Further, note that
  \begin{align*}
    \set{\bN\in\XXX(\TT)}{\ip{\dDiv\bN}v_\Omega=0 \quad \forall v\in P^1(\TT)} = \set{\bN\in\XXX(\TT)}{\dDiv\bN =0} % =: \mathbb{K}(\TT)
  \end{align*}
  since $\dDiv(\XXX(\TT))\subseteq P^1(\TT)$.
  That is, the discrete kernel is subspace of the continuous kernel, giving coercivity on the discrete kernel.
  The first two asserted error estimates then follow from the theory on mixed methods~\cite{BoffiBrezziFortin}. 

  The restricted quasi-optimality result for the error of the dual variable in weaker norm is also classical in mixed methods for second-order problems, see~\cite{BoffiBrezziFortin}. It can be seen as follows. Note that by~\eqref{eq:varform} and~\eqref{eq:mixed} we obtain
  \begin{align}\label{eq:galerkinorth}
    \ip{\Ctens^{-1}(\bM-\bM_\TT)}{\bQ}_\Omega - \ip{u-u_\TT}{\dDiv\bQ}_\Omega = 0 \quad\forall \bQ\in\XXX(\TT).
  \end{align}
  By~\eqref{eq:mixed:b} we have $\dDiv\bM_\TT = \Pi_\TT^1f =\Pi_\TT^1\dDiv\bM$. Let $\bN\in\XXX(\TT)$ be arbitrary with $\dDiv\bN = \Pi_\TT^1\dDiv\bM$. Setting $\bQ = \bM_\TT-\bN$ we see that $\dDiv\bQ = 0$ and infer that
  \begin{align*}
    0=\ip{\Ctens^{-1}(\bM-\bM_\TT)}{\bQ}_\Omega = \ip{\Ctens^{-1}(\bM-\bM_\TT)}{\bM_\TT-\bN}_\Omega.
  \end{align*}
  Using the latter identity we have that
  \begin{align*}
    \norm{\Ctens^{-1/2}(\bM-\bM_\TT)}\Omega^2 &= \ip{\Ctens^{-1}(\bM-\bM_\TT)}{\bM-\bM_\TT}_\Omega 
    \\
    &= \ip{\Ctens^{-1}(\bM-\bM_\TT)}{\bM-\bN}_\Omega \leq \norm{\Ctens^{-1/2}(\bM-\bM_\TT)}\Omega\norm{\Ctens^{-1/2}(\bM-\bN)}\Omega.
  \end{align*}
  This concludes the proof.
\end{proof}

\subsection{Accuracy enhancement by postprocessing}\label{sec:postproc}
In this section we postprocess solutions $(\bM_\TT,u_\TT)$ of mixed scheme~\eqref{eq:mixed} to achieve higher convergence rates in the primal variable under additional regularity assumptions. 
The postprocessing scheme is quite common and used in many works, see, e.g.~\cite{Stenberg91} for a similar technique or~\cite[Section~3.2]{ChenHuang23} in the context of mixed FEM for the biharmonic problem.
The result given below in Lemma~\ref{lem:superclose} is also found in similar form in~\cite[Remark~3.5]{ChenHuang23}.

Given $(\bM_\TT,u_\TT)\in \XXX(\TT)\times P^1(\TT)$, define $u_\TT^\star\in P^3(\TT)$ as the unique solution of 
\begin{subequations}\label{eq:postproc}
\begin{align*}
  \ip{\Ggrad u_\TT^\star}{\Ggrad v}_{\TT} &= \ip{\bM_\TT}{\Ggrad v}_{\TT} \quad\forall v\in P^3(\TT), \\
  \Pi_\TT^1 u_\TT^\star &= u_\TT.
\end{align*}
\end{subequations}

\begin{theorem}\label{thm:postproc}
  Suppose that $\Omega$ is convex and that the maximum interior angle at corner points of the domain $\Omega$ is smaller than $126.38^\circ$. Suppose that $\widehat\bu \in \trGgrad\Omega(H^4(\Omega))$, $f\in H^2(\TT)$ and that $\Ctens$ is the identity. Let $(\bM_\TT,u_\TT)$ denote the solution of~\eqref{eq:mixed} and define $u_\TT^\star\in P^3(\TT)$ by~\eqref{eq:postproc}. Then, 
  \begin{align*}
    \|u-u_\TT^\star\|_\Omega = \OO(\|h_\TT\|_{L^\infty(\Omega)}^4).
  \end{align*}
\end{theorem}

For the proof of Theorem~\ref{thm:postproc} we need the following super-closeness result. 
\begin{lemma}\label{lem:superclose}
  With the assumptions and notations from Theorem~\ref{thm:postproc} the estimate
  \begin{align*}
    \|\Pi_\TT^1-u_\TT\|_\Omega \lesssim \|h_\TT\|_{L^\infty(\Omega)}^2\|\bM-\bM_\TT\|_{\HdivDivset\Omega}
  \end{align*}
  holds.
\end{lemma}
\begin{proof}
  Let $v\in H_0^2(\Omega)$ denote the solution of 
  \begin{align*}
    \Delta^2 v = \Pi_\TT^1 u - u_\TT.
  \end{align*}
  By regularity results~\cite[Theorem~2]{BlumRannacher80} it follows that $v\in H^4(\Omega)$ with $\|v\|_{H^4(\Omega)}\lesssim \|\Pi_\TT^1u-u_\TT\|_\Omega$. Set $\bN = -\Ggrad v$. Then, 
  \begin{align*}
    \|\Pi_\TT^1 u-u_\TT\|^2 &= \ip{\Pi_\TT^1 u-u_\TT}{-\dDiv\bN}_\Omega + \ip{\bM-\bM_\TT}{\bN+\Ggrad v}_\Omega
    \\
    &= \ip{\bM-\bM_\TT}{\bN}_\Omega - \ip{u-u_\TT}{\dDiv\bN}_\Omega - \ip{\dDiv(\bM-\bM_\TT)}{-v}_\Omega
    \\
    &= \ip{\bM-\bM_\TT}{\bN-\bN_\TT}_\Omega - \ip{u-u_\TT}{\dDiv(\bN-\bN_\TT)}_\Omega + \ip{\dDiv(\bM-\bM_\TT)}{v}_\Omega
  \end{align*}
  for all $\bN_\TT\in \XXX(\TT)$.
  The last equality follows from Galerkin orthogonality. Choosing $\bN_\TT=\projHdDiv_\TT\bN$ and using the results from Proposition~\ref{prop:proj} we obtain with $\dDiv(\bM-\bM_\TT) = (1-\Pi_\TT^1)f$
the relations $\dDiv(\bN-\bN_\TT) = 0$ and
  \begin{align*}
    %\dDiv(\bN-\bN_\TT) = 0, \quad
    &|\ip{\dDiv(\bM-\bM_\TT)}{v}_\Omega| = |\ip{\dDiv(\bM-\bM_\TT)}{v-\Pi_\TT^1v}_\Omega| \\
    &\qquad\lesssim \|\dDiv(\bM-\bM_\TT)\|_\Omega \|h_\TT\|_{L^\infty(\Omega)}^2 \|\Pi_\TT^1 u-u_\TT\|_\Omega.
  \end{align*}
  Furthermore, 
  \begin{align*}
    |\ip{\bM-\bM_\TT}{\bN-\bN_\TT}_\Omega|&\leq \|\bM-\bM_\TT\|_\Omega \|h_\TT\|_{L^\infty(\Omega)}^2 \|\Pi_\TT^1 u-u_\TT\|_\Omega.
  \end{align*}
  Putting all the estimates together concludes the proof. 
\end{proof}

\begin{proof}[Proof of Theorem~\ref{thm:postproc}]
  The proof is quite standard for the proposed postprocessing scheme, see, e.g.,~\cite{Stenberg91}. 
  Therefore, we only give some details. Note that by the assumptions of the theorem,
  %it follows that
  $u \in H^4(\Omega)$, so that
  \begin{align*}
    \|\bM-\bM_\TT\|_{\HdivDivset\Omega} = \OO(\|h_\TT\|_{L^\infty(\Omega)}^2).
  \end{align*}
  Then, together with Lemma~\ref{lem:superclose} we conclude
  \begin{align*}
    \|u-u_\TT^\star\|_\Omega &\leq \|\Pi_\TT^1(u-u_\TT^\star)\|_\Omega + \|(1-\Pi_\TT^1)(u-u_\TT^\star)\|_\Omega
    \\
    &= \|\Pi_\TT^1u-u_\TT\|_\Omega + \|(1-\Pi_\TT^1)(u-u_\TT^\star)\|_\Omega
    \\
    &\lesssim \|h_\TT\|_{L^\infty(\Omega)}^4 + \|h_\TT\|_{L^\infty(\Omega)}^2\|\Ggrad(u-u_\TT^\star)\|_\Omega.
  \end{align*}
  Let $w_\TT\in P^3(\TT)$ satisfy $\ip{\Ggrad w_\TT}{\Ggrad v}_\TT = \ip{\bM}{\Ggrad v}_\TT$ for all $v\in P^3(\TT)$. Since $\Ggrad u = \bM$ we infer that $\|\Ggrad(u-w_\TT)\|_\Omega \lesssim \|h_\TT\|_{L^\infty(\Omega)}^2\|u\|_{H^4(\Omega)}$. This gives
  \begin{align*}
    \|\Ggrad(u-u_\TT^\star)\|_\Omega &\leq \|\Ggrad(u-w_\TT)\|_\Omega + \|\Ggrad(w_\TT-u_\TT^\star)\|_\Omega
    \\
    &\lesssim  \|h_\TT\|_{L^\infty(\Omega)}^2\|u\|_{H^4(\Omega)} + \|\bM-\bM_\TT\|_\Omega.
  \end{align*}
  Combining all estimates finishes the proof.
\end{proof}

\section{A posteriori error estimation}\label{sec:aposteriori}
In this section we derive an error estimator for the bending moments in the $L_2(\Omega)$-norm.
We first define the estimator and then state the main result on its reliability and efficiency. Proofs are postponed to Sections~\ref{sec:proofRel} and~\ref{sec:proofEff}.
While the following analysis is mostly independent of space $\XXX(\TT)$, we assume
for simplicity that $\TT=\TT_\triangle$ only consists of triangles.
Moreover, we assume that $\Ctens$ is the identity.

For the remainder of this section, let $u\in H^2(\Omega)$ denote the solution of~\eqref{eq:biharm}, so that $(\bM,u) = (\Ggrad u,u)$ is the solution of variational formulation~\eqref{eq:varform}. 
We consider the solution component $\bM_\TT$ of the mixed FEM~\eqref{eq:mixed}
and denote by $\jump{\bM_\TT\tangential}|_E$ the jump of $\bM_\TT\tangential$ over an interior edge
$E\in\edges_\Omega$. We make the standing assumption that $\bg:=\bM\tangential|_\Gamma = (\Ggrad u)\tangential|_\Gamma\in L_2(\Gamma)$. The local error indicators are given by
\begin{align*}
  \osc_K &:= h_K^2\norm{(1-\Pi_K^1)f}K, \\
  \nu_{\TT}(K)^2 &:= h_K^2 \norm{\Rot\bM_\TT}K^2 + h_K\norm{(1-\Pi_{\edges_K}^0)\jump{\bM_\TT\tangential}}{\partial K\setminus\Gamma}^2 
  \\
  &\qquad + h_K\norm{(1-\Pi_{\edges_K}^0)(\bM_\TT\tangential-\bg)}{\partial K\cap \Gamma}^2 + \osc_K^2.
\end{align*}
Here, $\Pi_{\edges_K}^0 v|_E := \Pi_E^0 v$ denotes the $L_2(E)$-orthogonal projection onto $P^0(E)$.
The (squared) total estimator is defined as the sum of the (squared) local contributions, i.e., 
\begin{align*}
  \nu_\TT^2 := \sum_{K\in\TT} \nu_\TT(K)^2.
\end{align*}
A similar error estimator has been derived in~\cite{HuangHuangXu11} for the Hellan--Herrmann--Johnson method with homogeneous clamped boundary conditions.
\begin{theorem}\label{thm:aposteriori}
  Let $\Omega$ denote a simply connected domain with triangulation $\TT$. 
  Under the aforegoing assumptions the estimator $\nu_\TT$ is reliable, i.e., 
  \begin{align*}
    \norm{\bM-\bM_\TT}{\Omega} \lesssim \nu_\TT.
  \end{align*}
  Let $K\in\TT$. If $\bg=\bM\tangential|_\Gamma$ is a piecewise polynomial of degree $\leq k$ on $\edges_K\cap \edges_\Gamma$, then
  \begin{align*}
    \nu_\TT(K) \lesssim \norm{\bM-\bM_\TT}{\Patch_\TT(K)} + \osc_{K},
  \end{align*}
  where the involved constant depends on $k$, but not on the particular $\bg$.
\end{theorem}

\subsection{Some tools}
We follow the general ideas from Carstensen~\cite{CC97} for proving Theorem~\ref{thm:aposteriori},
see also~\cite{HuangHuangXu11} for a similar a posteriori analysis for the Hellan--Herrmann--Johnson method.
One of the main ingredients in the proofs is a Helmholtz decomposition of vector fields. 
In this work we consider the following Helmholtz-type decomposition of tensor fields. 
It is based on~\cite[Section~4]{RafetsederZulehner18}, see also~\cite[Section~2.4]{ChenHuang18} for a similar result with vanishing boundary conditions.
\begin{lemma}\label{lem:helmholtz}
  Let $v\in H^2(\Omega)$ be given and let $\Omega$ be simply connected. 
  For any $\bN\in\HdivDivset\Omega$ there exists $\bq\in\bH^1(\Omega)_{/\RT^0(\Omega)}$ such that
  \begin{align*}
    \bN = \Ggrad p + \symCurl \bq
  \end{align*}
  with $p\in H^2(\Omega)$ denoting the unique solution of 
  \begin{align*}
    \dDiv \Ggrad p &= \dDiv\bN \quad\text{in }\Omega, \\
    \trGgrad{\Omega}p &= \trGgrad{\Omega}v.
  \end{align*}
  In particular, $\norm{\bq}{\bH^1(\Omega)}\lesssim \norm{\symCurl\bq}{\Omega}$.
\end{lemma}
\begin{proof}
  Defining $p$ as in the statement we have that $\dDiv(\bN-\Ggrad p)=0$. By~\cite[Lemma~4.1]{RafetsederZulehner18} there exists $\bq\in\bH^1(\Omega)$ with 
  \begin{align*}
    \symCurl\bq = \bN-\Ggrad p.
  \end{align*}
  Function $\bq$ is unique up to an element in $\RT^0(\Omega)$. 
  The last estimate is a Korn inequality. 
  This can be seen with the arguments from~\cite[Section~4]{RafetsederZulehner18}, see particularly~\cite[Remark~4.3]{RafetsederZulehner18}. We can thus choose $\bq\in\bH^1(\Omega)_{/\RT^0(\Omega)}$. 
  This finishes the proof.
\end{proof}
For the remainder of this section we use the decomposition from Lemma~\ref{lem:helmholtz} with $v=u$
and $\bN = \bM_\TT$, giving
\begin{align}
  \bM_\TT &= \Ggrad p + \symCurl\bq. \label{eq:decompositions:MT}
\end{align}
Integration by parts, $\dDiv\symCurl\bq =0$, and $u-p\in H_0^2(\Omega)$ show that
\begin{align}\label{eq:ortoUP}
  \ip{\Ggrad(u-p)}{\symCurl\bq}_\Omega &= 0.
\end{align}
Recalling that $\bM=\Ggrad u$ we directly obtain 
\begin{align}\label{eq:splittingError}
  \norm{\bM-\bM_\TT}{\Omega}^2 = \norm{\Ggrad(u-p)}{\Omega}^2 + \norm{\symCurl\bq}{\Omega}^2.
\end{align}
Let $\projSZtilde_\TT\colon H^1(\Omega)\to P^2(\TT)\cap H^1(\Omega)$ denote a quasi-interpolator with
\begin{align*}
  \norm{\nabla\projSZtilde_\TT v}K\lesssim \norm{\nabla v}{\Patch_\TT(K)}, 
  \quad
  \norm{(1-\projSZtilde_\TT v)}K\lesssim h_K \norm{\nabla v}{\Patch_\TT(K)}
\end{align*}
for all $K\in\TT$, $v\in H^1(\Omega)$. The trace inequality then shows that
\begin{align*}
  \norm{(1-\projSZtilde_\TT)v}E \lesssim h_E^{1/2} \norm{\nabla v}{\Patch_\TT(K_E)}
\end{align*}
where $K_E\in\TT$ is any element with $E\in\edges_K$.
An example for such an operator is the Scott--Zhang operator~\cite{SZ_90} which is also a projection, or Cl\'ement's operator~\cite{Clement75}.
We consider a slight modification of this operator where we add correction terms to ensure orthogonality on the edges. To that end let $\eta_E \in P^2(\TT)\cap H^1(\Omega)$ denote the edge bubble given as the product of the barycentric coordinate functions of the two vertices of edge $E$.
\begin{lemma}
  Consider the operator $\projSZ_\TT\colon H^1(\Omega)\to P^2(\TT)\cap H^1(\Omega)$, 
  \begin{align*}
    \projSZ_\TT v = \projSZtilde_\TT v + \sum_{E\in\edges} \frac{\dual{(1-\projSZtilde_\TT)v}1_E}{\dual{\eta_E}{1}_E}\eta_E.
  \end{align*}
  It satisfies
  \begin{align*}
    \norm{\nabla\projSZ_\TT v}K\lesssim \norm{\nabla v}{\Patch_\TT(K)}, 
    \quad
    \norm{(1-\projSZ_\TT )v}K\lesssim h_K \norm{\nabla v}{\Patch_\TT(K)}
    \quad\forall K\in\TT
  \end{align*}
  and
  \begin{align*}
    \norm{(1-\projSZ_\TT)v}E \lesssim h_E^{1/2}\norm{\nabla v}{\Patch_\TT(K_E)}, \quad
    \dual{(1-\projSZ_\TT)v}{1}_E = 0 \quad\forall E\in\edges.
  \end{align*}
  Here, for any $E\in\edges$, $K_E\in\TT$ is an element with $E\in\edges_K$.
\end{lemma}
\begin{proof}
  The orthogonality relation follows directly by definition of the operator. Let $v\in H^1(\Omega)$ and $E\in\edges$ be given. Then, relation
  \begin{align*}
    \dual{(1-\projSZ_\TT)v}{1}_E = \dual{(1-\projSZtilde_\TT)v}1_E - \frac{\dual{(1-\projSZtilde_\TT)v}1_E}{\dual{\eta_E}1_E} \dual{\eta_E}1_E =0
  \end{align*}
  holds.
  The other properties follow from the ones of $\projSZtilde_\TT$ and scaling arguments. Let $K\in\TT$ and $v\in H^1(\Omega)$ be given. First, 
  \begin{align*}
    \norm{(1-\projSZ_\TT)v}K &\leq \norm{(1-\projSZtilde_\TT)v}K + \sum_{E\in\edges_K} \frac{\norm{(1-\projSZtilde_\TT)v}E\norm{1}E}{|\dual{\eta_E}1_E|}\norm{\eta_E}K 
    \\ &\lesssim \norm{(1-\projSZtilde_\TT)v}K + \sum_{E\in\edges_K} h_E^{1/2} \norm{(1-\projSZtilde_\TT)v}E
    \lesssim h_K \norm{\nabla v}{\Patch_\TT(K)}.
  \end{align*}
  Second, 
  \begin{align*}
    \norm{\nabla\projSZ_\TT v}K & \lesssim \norm{\nabla\projSZtilde_\TT v}K + \sum_{E\in\edges_K} h_E^{-1/2}\norm{(1-\projSZtilde_\TT)v}E \lesssim \norm{\nabla v}{\Patch_\TT(K)}.
  \end{align*}
  Finally, 
  \begin{align*}
    \norm{(1-\projSZ_\TT)v}E \lesssim \norm{(1-\projSZtilde_\TT)v}E + \norm{(1-\projSZtilde_\TT)v}E \lesssim 
    h_E^{1/2}\norm{\nabla v}{\Patch_\TT(K)}
  \end{align*}
  holds for any $E\in\edges_K$, which concludes the proof.
\end{proof}

\subsection{Proof of reliability in Theorem~\ref{thm:aposteriori}}\label{sec:proofRel}
For the proof of reliability we start with~\eqref{eq:splittingError}.
From~\eqref{eq:mixed:b} and $\dDiv\bM_\TT\in P^1(\TT)$ we find that $\dDiv\bM_\TT = \Pi_\TT^1 f$. 
For the first term on the right-hand side of~\eqref{eq:splittingError} we therefore get by integration by parts
\begin{align*}
  \norm{\Ggrad(u-p)}\Omega^2 &= \ip{\Ggrad(u-p)}{\Ggrad(u-p)}_\Omega = \ip{u-p}{\dDiv(\bM-\bM_\TT)}_\Omega 
  \\
  &= \ip{u-p}{(1-\Pi_\TT^1)f}_\Omega 
  = \ip{(1-\Pi_\TT^1)(u-p)}{(1-\Pi_\TT^1)f}_\Omega 
  \\
  &\lesssim \norm{\Ggrad(u-p)}\Omega\norm{h_\TT^2(1-\Pi_\TT^1)f}\Omega.
\end{align*}
This proves that $\norm{\Ggrad(u-p)}\Omega \lesssim \norm{h_\TT^2(1-\Pi_\TT^1)f}\Omega$.

For the second term on the right-hand side of~\eqref{eq:splittingError} let $\bq_\TT\in\bP^2(\TT)\cap\bH^1(\Omega)$ be arbitrary. Note that $\symCurl\bq_\TT\in \XXX(\TT)$. 
Further note that $\ip{\Ggrad(u-p)}{\symCurl\bq_\TT}_\Omega = 0$ since $u-p\in H_0^2(\Omega)$ and $\dDiv\symCurl\bq_\TT=0$.
Using Galerkin orthogonality~\eqref{eq:galerkinorth} we find that
\begin{align*}
  \ip{-\symCurl\bq}{\symCurl\bq_\TT}_\Omega = \ip{\Ggrad(u-p)-\symCurl\bq}{\symCurl\bq_\TT}_\Omega = \ip{\bM-\bM_\TT}{\symCurl\bq_\TT}_\Omega= 0.
\end{align*}
Then, using $\ip{\Ggrad(u-p)}{\symCurl(\bq-\bq_\TT)}_\Omega = 0$ and integrating by parts,
\begin{align*}
  \norm{\symCurl\bq}{\Omega}^2 &= \ip{\symCurl\bq}{\symCurl(\bq-\bq_\TT)}_\Omega
  \\ &= \ip{\Ggrad p+\symCurl\bq}{\symCurl(\bq-\bq_\TT)}_\Omega - \ip{\Ggrad p}{\symCurl(\bq-\bq_\TT)}_\Omega
  \\
  &= \ip{\bM_\TT}{\symCurl(\bq-\bq_\TT)}_\Omega - \ip{\Ggrad u}{\symCurl(\bq-\bq_\TT)}_\Omega
  \\
  &= \ip{\bM_\TT}{\symCurl(\bq-\bq_\TT)}_\Omega - \ip{\Ggrad u}{\Curl(\bq-\bq_\TT)}_\Omega
  \\
  &= \ip{\bM_\TT}{\symCurl(\bq-\bq_\TT)}_\Omega + \dual{\bg}{\bq-\bq_\TT}_\Gamma.
\end{align*}
We rewrite the first term by first noting that 
\begin{align*}
  \ip{\bM_\TT}{\symCurl(\bq-\bq_\TT)}_K &= \ip{\bM_\TT}{\Curl(\bq-\bq_\TT)}_K =
  \ip{\Rot\bM_\TT}{\bq-\bq_\TT}_K - \dual{\bM_\TT\tangential}{\bq-\bq_\TT}_{\partial K}.
\end{align*}
Summing over all elements and combining terms on interior edges we find that
\begin{align*}
  \ip{\bM_\TT}{\symCurl(\bq-\bq_\TT)}_\Omega + \dual{\bg}{\bq-\bq_\TT}_\Gamma
  &= \sum_{K\in\TT} \ip{\Rot\bM_\TT}{\bq-\bq_\TT}_K + \sum_{E\in\edges_\Omega} \dual{\jump{\bM_\TT\tangential}}{\bq-\bq_\TT}_E 
  \\ 
  &\qquad + \sum_{E\in\edges_\Gamma} \dual{\bg-\bM_\TT\tangential}{\bq-\bq_\TT}_E.
\end{align*}
Choosing $\bq_\TT = \projSZ_\TT\bq$ and using the properties of $\projSZ_\TT$ as well as $\norm{\bq}{\bH^1(\Omega)}\lesssim \norm{\symCurl\bq}\Omega$ we obtain
\begin{align*}
  \sum_{K\in\TT} \ip{\Rot\bM_\TT}{\bq-\bq_\TT}_K \lesssim \norm{h_\TT\Rot_\TT\bM_\TT}{\Omega}\norm{\symCurl\bq}\Omega
\end{align*}
and
\begin{align*}
  \sum_{E\in\edges_\Omega} \dual{\jump{\bM_\TT\tangential}}{\bq-\bq_\TT}_E &
  = \sum_{E\in\edges_\Omega} \dual{(1-\Pi_E^0)\jump{\bM_\TT\tangential}}{\bq-\bq_\TT}_E
  \\ &\lesssim \sqrt{\sum_{E\in\edges_\Omega} h_E\norm{(1-\Pi_E^0)\jump{\bM_\TT\tangential}}E^2}\norm{\symCurl\bq}\Omega
  \\
  &\lesssim \sqrt{\sum_{K\in\TT}h_K\norm{(1-\Pi_{\edges_K}^0)\jump{\bM_\TT\tangential}}{\partial K\setminus\Gamma}^2}\norm{\symCurl\bq}\Omega.
\end{align*}
For the remaining boundary terms similar arguments show that
\begin{align*}
  \sum_{E\in\edges_\Gamma} \dual{\bg-\bM_\TT\tangential}{\bq-\bq_\TT}_E
  &\lesssim \sqrt{\sum_{E\in\edges_\Gamma} h_E\norm{(1-\Pi_E^0)(\bg-\bM_\TT\tangential)}E^2} \norm{\symCurl\bq}\Omega
  \\
  &\lesssim \sqrt{\sum_{K\in \TT} h_K\norm{(1-\Pi_{\edges_K}^0)(\bg-\bM_\TT\tangential)}{\partial K\cap\Gamma}^2} \norm{\symCurl\bq}\Omega.
\end{align*}
Putting all the estimates together we conclude that
\begin{align*}
  \norm{\symCurl\bq}\Omega^2 &\lesssim \norm{h_\TT\Rot_\TT\bM_\TT}{\Omega}^2 + \sum_{K\in\TT}h_K\norm{(1-\Pi_{\edges_K}^0)\jump{\bM_\TT\tangential}}{\partial K\setminus\Gamma}^2 \\
  &\qquad + \sum_{K\in \TT} h_K\norm{(1-\Pi_{\edges_K}^0)(\bg-\bM_\TT\tangential)}{\partial K\cap\Gamma}^2.
\end{align*}
Together with the estimate $\norm{\Ggrad(u-p)}\Omega \lesssim \norm{h_\TT^2(1-\Pi_\TT^1)f}\Omega$ established before, this finishes the proof. 

\subsection{Proof of efficiency in Theorem~\ref{thm:aposteriori}}\label{sec:proofEff}
Local efficiency is shown by using Verf\"{u}rth's bubble function technique~\cite{Verfuerth94}. For an application of this technique to mixed FEM for scalar second-order elliptic equations we refer to~\cite{CC97}. 

We divide the proof into three steps presented in the next three lemmas. Combining these results together with the simple estimate
\begin{align*}
  \norm{(1-\Pi_E^0)g}E \leq \norm{g}E \quad\forall g\in L_2(E), \, E\in\edges
\end{align*}
proves the efficiency bound from Theorem~\ref{thm:aposteriori}. 

\begin{lemma}\label{lem:eff:vol}
  The estimate
  \begin{align*}
    h_K\norm{\Rot\bM_\TT}K \lesssim \norm{\bM-\bM_\TT}K
  \end{align*}
  holds for all $K\in\TT$.
\end{lemma}
\begin{proof}
  Let $\eta_K$ denote the element bubble function, i.e., the product of the barycentric coordinate functions. Norm equivalence in finite-dimensional spaces and scaling arguments, $\Rot\bM = \Rot\Ggrad u = 0$, and integration by parts prove that
  \begin{align*}
    \norm{\Rot\bM_\TT}K^2 &\eqsim \ip{\Rot\bM_\TT}{\eta_K\Rot\bM_\TT}_K = \ip{\Rot(\bM_\TT-\bM)}{\eta_K\Rot\bM_\TT}_K
    \\
    &= \ip{\bM_\TT-\bM}{\Curl(\eta_K\Rot\bM_\TT)}_K \leq \norm{\bM_\TT-\bM}K\norm{\Curl(\eta_K\Rot\bM_\TT)}K.
  \end{align*}
  Noting that $\eta_K\Rot\bM_\TT|_K$ is a polynomial, we use the inverse inequality $\norm{\Curl(\eta_K\Rot\bM_\TT)}K\lesssim h_K^{-1}\norm{\eta_K\Rot\bM_\TT}K$ to see that
  \begin{align*}
    \norm{\Rot\bM_\TT}K^2 \lesssim \norm{\bM_\TT-\bM}K h_K^{-1}\norm{\eta_K\Rot\bM_\TT}K \leq h_K^{-1} \norm{\bM_\TT-\bM}K \norm{\Rot\bM_\TT}K.
  \end{align*}
  Dividing by $\norm{\Rot\bM_\TT}K$ and multiplying with $h_K$ finishes the proof.
\end{proof}

\begin{lemma}\label{lem:eff:jump}
  The estimate
  \begin{align*}
    h_K^{1/2}\norm{\jump{\bM_\TT\tangential}}{\partial K\setminus \Gamma} \lesssim \norm{\bM-\bM_\TT}{\Patch_\TT(K)}
  \end{align*}
  holds for all $K\in\TT$.
\end{lemma}
\begin{proof}
  Let $E\in\edges_K\cap \edges_\Omega$ denote some interior edge of $K$. There exists a unique $K'\in\TT$, $K'\neq K$, such that $\overline K\cap \overline K' = \overline E$. Denote by $\patch_E = \{K,K'\}$ and $\Patch_E$ the associated domain. There exists an extension operator $P\colon C(E)\to C(\Patch_E)$ (see~\cite[Proof of Lemma~6.2]{CC97} and references therein) such that $P\sigma|_E = \sigma$ for polynomials of degree $\leq k$, and
  \begin{align*}
    h_E \norm{\sigma}{E}^2 \eqsim h_E \norm{\eta_E^{1/2}\sigma}{E}^2 \eqsim  \norm{\eta_E^{1/2}P\sigma}{\Patch_E}^2, \qquad
    \norm{\nabla(\eta_EP\sigma)}{\Patch_E} \lesssim h_E^{-1/2}\norm{\sigma}{E}.
  \end{align*}
  We apply this operator to each component of $\bsigma = \jump{\bM_\TT\tangential}|_E$. 
  We use notation $\norm{\cdot}{\patch_E}^2 =\sum_{K\in\patch_E} \norm{\cdot}{K}^2$ for broken norms and a similar one for the corresponding inner product.
  Using $\jump{\bM\tangential}|_E = 0$, integration by parts, $\Rot\bM = 0$, and the properties of the above mentioned extension operator we find that
  \begin{align*}
    \norm{\jump{\bM_\TT\tangential}}E^2 &\eqsim \dual{\jump{\bM_\TT\tangential}}{\eta_EP\bsigma}_E 
    = \dual{\jump{(\bM_\TT-\bM)\tangential}}{\eta_EP\bsigma}_E
    \\
    &= \ip{\Rot(\bM_\TT-\bM)}{\eta_EP\bsigma}_{\patch_E} - \ip{\bM_\TT-\bM}{\Curl(\eta_EP\bsigma)}_{\Patch_E}
    \\
    &\lesssim \norm{\eta_E^{1/2}\Rot(\bM_\TT-\bM)}{\patch_E}\norm{\eta_E^{1/2}P\bsigma}{\Patch_E} + \norm{\bM_\TT-\bM}{\Patch_E} h_E^{-1/2}\norm{\bsigma}{E}.
    \\
    &\lesssim \norm{\Rot\bM_\TT}{\patch_E} h_E^{1/2}\norm{\bsigma}E + \norm{\bM_\TT-\bM}{\Patch_E} h_E^{-1/2}\norm{\bsigma}E.
  \end{align*}
  Then, $h_E\eqsim h_K$ gives
  \begin{align*}
    h_K^{1/2}\norm{\jump{\bM_\TT\tangential}}E \lesssim h_K \norm{\Rot\bM_\TT}{\patch_E} + \norm{\bM_\TT-\bM}{\Patch_E}.
  \end{align*}
  Application of Lemma~\ref{lem:eff:vol} to bound $h_K \norm{\Rot\bM_\TT}{\patch_E}$ and summation over all interior edges of $K$ finishes the proof.
\end{proof}

\begin{lemma}\label{lem:eff:bou}
  Assuming that $\bg|_E \in P^k(E)$ for all $E\in\edges_K\cap \edges_\Gamma$, estimate
  \begin{align*}
    h_K^{1/2}\norm{\bg-\bM_\TT\tangential}{\partial K\cap \Gamma} \lesssim \norm{\bM-\bM_\TT}{\Patch_\TT(K)}
    %+ h_K^{1/2}\norm{(1-\Pi_{\edges_K}^3)\partial_{\tangential}\bg}{\partial K\cap \Gamma}.
  \end{align*}
  holds for all $K\in\TT$.
  The involved constant depends on the polynomial degree $k$ but is independent of $\bg$.
\end{lemma}
\begin{proof}
  Noting that $(\bg-\bM_\TT\tangential)|_E = (\bM-\bM_\TT)\tangential|_E$ is a polynomial on all edges $E\in \edges_K\cap\edges_\Gamma$ we may argue as before in Lemma~\ref{lem:eff:jump}, and omit further details.
\end{proof}

%%%
\section{Numerical experiments}\label{sec:ex}
In this section we present two numerical experiments for the mixed scheme~\eqref{eq:mixed} with $\Ctens$ being the identity. The first one in Section~\ref{sec:ex:smooth} considers a smooth solution in a convex domain and the second one in Section~\ref{sec:ex:sing} a typical singularity solution in a non-convex domain.
We consider sequences of uniformly refined meshes, $\TT_0,\TT_1,\dots$, where $\TT_{\ell+1}$ is constructed from $\TT_\ell$ by bisecting each triangle twice according to the newest vertex bisection rule (NVB) if $\TT_0$ is a mesh of triangles. In the case that $\TT_0$ is a mesh of parallelograms each element is divided into four by connecting the midpoints of opposite edges. 
Note that $\norm{h_{\TT_\ell}}{L_\infty(\Omega)} \eqsim (\#\TT_\ell)^{-1/2}$.
For the second experiment we additionally consider a sequence of locally refined meshes where we also use NVB for refining elements, but mark elements for refining according to the following simple adaptive loop:
\\
\textbf{Input:} Initial triangulation $\TT_0$, data $\widehat\bu$, $f$, marking parameter $\theta\in(0,1)$, and counter $\ell:=0$.
\\
Repeat the following steps:
\begin{itemize}
  \item \textbf{Solve:} compute solution to mixed FEM~\eqref{eq:mixed} on mesh $\TT_\ell$.
  \item \textbf{Estimate:} compute local error indicators $\nu_{\TT_\ell}(K)$ for all $K\in\TT_\ell$.
  \item \textbf{Mark:} mark elements for refinement according to the bulk criterion: Find a minimal set $\marked_\ell\subseteq \TT_\ell$ such that
    \begin{align*}
      \theta \nu_{\TT_\ell}^2 \leq \sum_{K\in\marked_\ell} \nu_{\TT_\ell}(K)^2.
    \end{align*}
  \item \textbf{Refine:} refine mesh $\TT_\ell$ to obtain $\TT_{\ell+1}$ such that at least all marked elements are refined and update counter $\ell\mapsto \ell+1$.
\end{itemize}
\textbf{Output:} Sequence of meshes $(\TT_\ell)_\ell$ and solutions $\big( (\bM_\ell,u_\ell)\big)_\ell$.

\subsection{Convex domain}\label{sec:ex:smooth}
\begin{figure}
  \begin{center}
    \input{ExampleSmooth}
    \caption{Errors and estimator for the smooth solution from Section~\ref{sec:ex:smooth};
             uniformly refined meshes with triangles (left) and squares (right).
             The black dotted line indicates $\OO\big((\#\TT)^{-1}\big)$.}
    \label{fig:smooth}
  \end{center}
\end{figure}
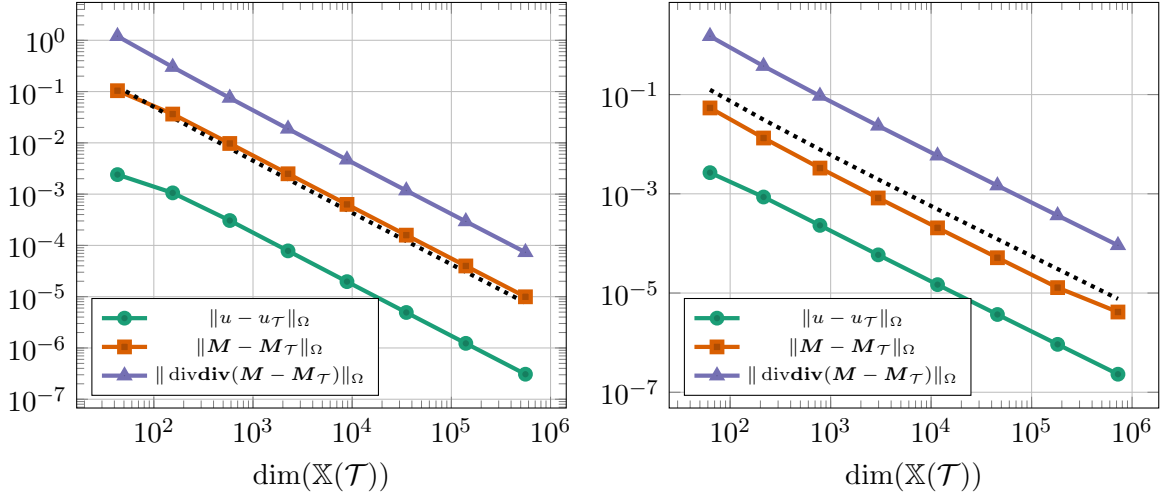

We consider the manufactured solution $u(x,y) = x^2y^2(1-x)(1-y)$ which satisfies the biharmonic problem~\eqref{eq:biharm} in $\Omega=(0,1)^2$ with $f = \Delta^2 u$ and $\widehat\bu = \trGgrad\Omega u$. Note that $\widehat\bu\neq 0$ because
\begin{align*}
 \partial_\normal u|_\Gamma(x,y) = 
  \begin{cases}
    -y^2(1-y) & x = 1, \\
    -x^2(1-x) & y =1, \\
    0 & \text{else}.
  \end{cases}
\end{align*}
We consider two sequences of uniform refined meshes. The first sequence uses the initial triangulation of $\Omega$ into four triangles $T_j = \conv\{z_j,z_{j+1},z_5\}$, $j=1,2,3,4$, and $z_1 = (0,0)$, $z_2=(1,0)$, $z_3 = (1,1)$, $z_4=(0,1)$, $z_5 = (\tfrac12,\tfrac12)$.
The second sequence uses the initial triangulation of $\Omega$ into four squares $T_j = \conv\{z_j,(z_j+z_{j+1})/2,z_5,(z_{j-1}+z_j)/2\}$, $j=1,2,3,4$, with $z_0 = z_4$ and $z_5 = z_1$.
Given that solution $u$ is smooth one combines Theorem~\ref{thm:mixed} and Proposition~\ref{prop:proj} to see that
\begin{align*}
  \norm{\bM-\bM_\TT}{\HdivDivset\Omega} + \norm{u-u_\TT}\Omega \lesssim \norm{h_\TT}{L_\infty(\Omega)}^2 = \OO\big( (\#\TT)^{-1} \big).
\end{align*}
Figure~\ref{fig:smooth} shows that these rates are indeed observed in the experiments. 
In particular, we find that all the error quantities $\norm{\bM-\bM_\TT}\Omega$, $\norm{u-u_\TT}\Omega$, and $\norm{\dDiv(\bM-\bM_\TT)}\Omega$ converge at the predicted rate for triangular as well as parallelogram grids.
% In addition, the error estimator $\nu_{\TT}$ is visualized in the same figure. Note that the oscillation term is of higher order, i.e., 
% \begin{align*}
%   \osc_\TT = \norm{h_\TT^2(1-\Pi_\TT^1)f}\Omega \lesssim \norm{h_\TT}{L_\infty(\Omega)}^4 = \OO\big( (\#\TT)^{-2} \big). 
% \end{align*}

\subsection{Non-convex domain}\label{sec:ex:sing}
\begin{figure}
  \begin{center}
    \input{ExampleSingInitMesh}
    \caption{Domain and initial meshes for the problem from Section~\ref{sec:ex:sing}.}
    \label{fig:sing:initMesh}
  \end{center}
\end{figure}
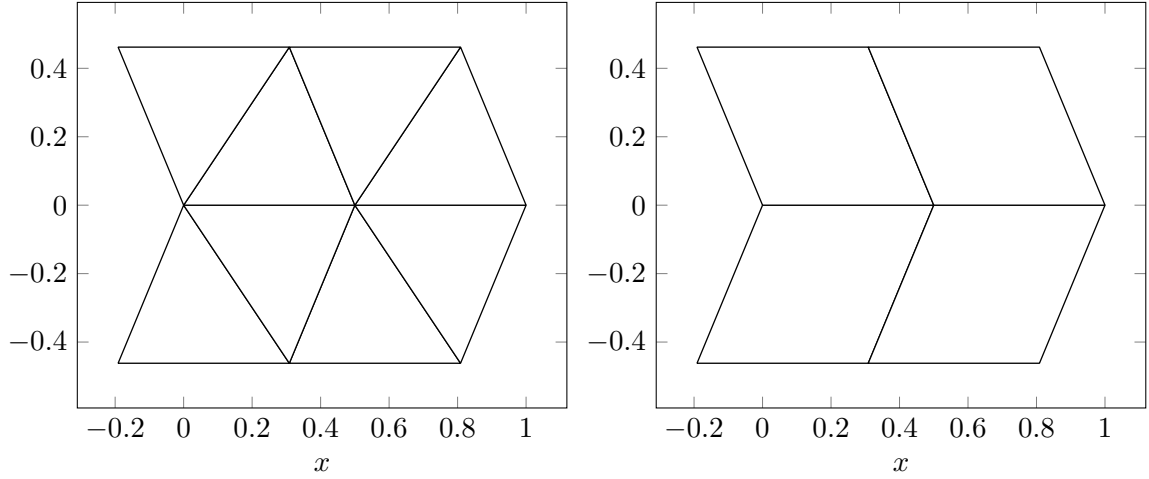

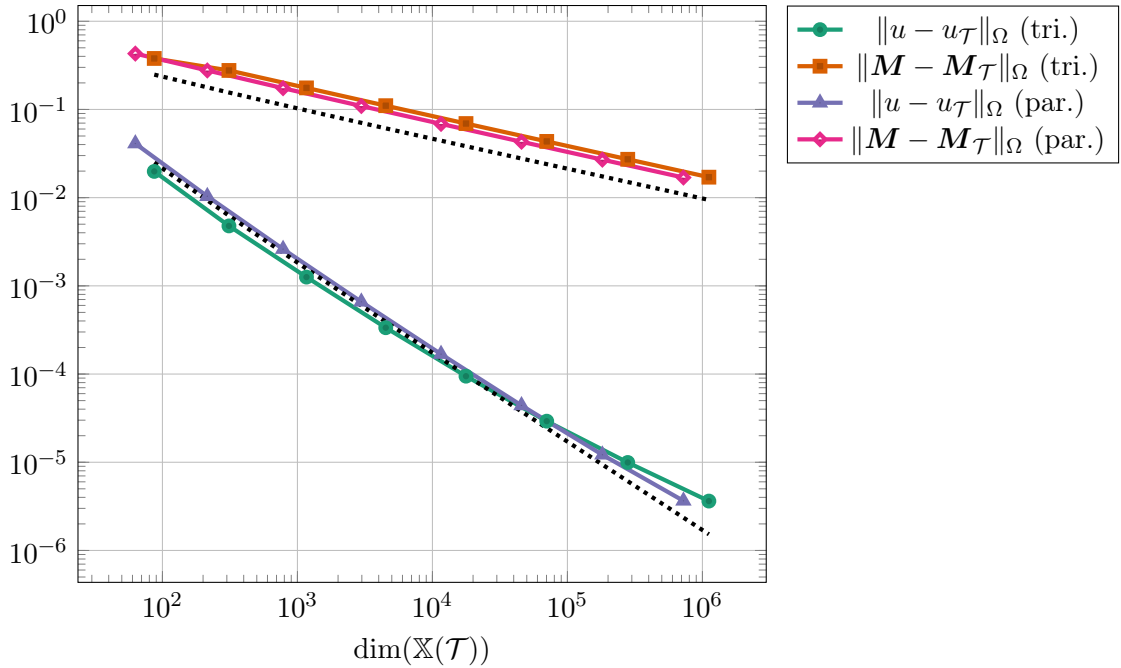
\begin{figure}
  \begin{center}
    \input{ExampleSing}
    \caption{Errors for the singular problem from Section~\ref{sec:ex:sing}
              on sequences of uniformly refined triangular and parallelogram meshes.
              The two black dotted lines indicate $\OO\big((\#\TT)^{-1}\big)$
              and $\OO\big((\#\TT)^{-0.337}\big)$.}
    \label{fig:sing}
  \end{center}
\end{figure}

For this experiment we consider the manufactured solution from~\cite[Section~6.2.2]{KLove1} with
domain $\Omega$ as given in Figure~\ref{fig:sing:initMesh}.
The initial triangulation with triangles is shown on the left and with parallelograms on the right.
In comparison to~\cite[Section~6.2.2]{KLove1} we use a slightly modified domain that can be
decomposed into triangles as well as parallelograms.
The domain has a re-entrant corner at the origin with interior angle $\tfrac54\pi$.
The manufactured solution is given by
\begin{align*}
  u(x,y) = r^{1+s}\big( \cos( (1+s)\varphi) + C\cos( (1-s)\varphi)\big)
\end{align*}
where $(r,\varphi)$ denote polar coordinates with $\varphi\in(-\pi,\pi]$.
We choose $s$ and $C$ such that $u\in H^{2+s-\epsilon}(\Omega)$ ($\epsilon>0$) is a typical singularity function of the biharmonic problem with vanishing traces $u$ and $\nabla u$ on the two edges that meet at the origin. 
For the present case this requires to set $s\approx 0.673583$ and $C\approx 1.23459$. The singularity function satisfies~\eqref{eq:biharm} with $f=0$. 
Furthermore, note that $\widehat\bu=\trGgrad\Omega u \neq 0$, but $\bg = \bM\tangential$ is smooth on all boundary parts. 
Additionally, we stress that $\bM\in \HHH^s(\Omega)$, but $\Div\bM\neq L_2(\Omega)$. In view of a priori approximation results we thus expect the reduced convergence
\begin{align*}
  \norm{u-u_\TT}\Omega + \norm{\bM-\bM_\TT}{\Omega} = \OO\big( (\#\TT)^{-s/2} \big)
\end{align*}
on a sequence of uniformly refined meshes. This is indeed observed in Figure~\ref{fig:sing} for triangular as well as parallelogram grids.

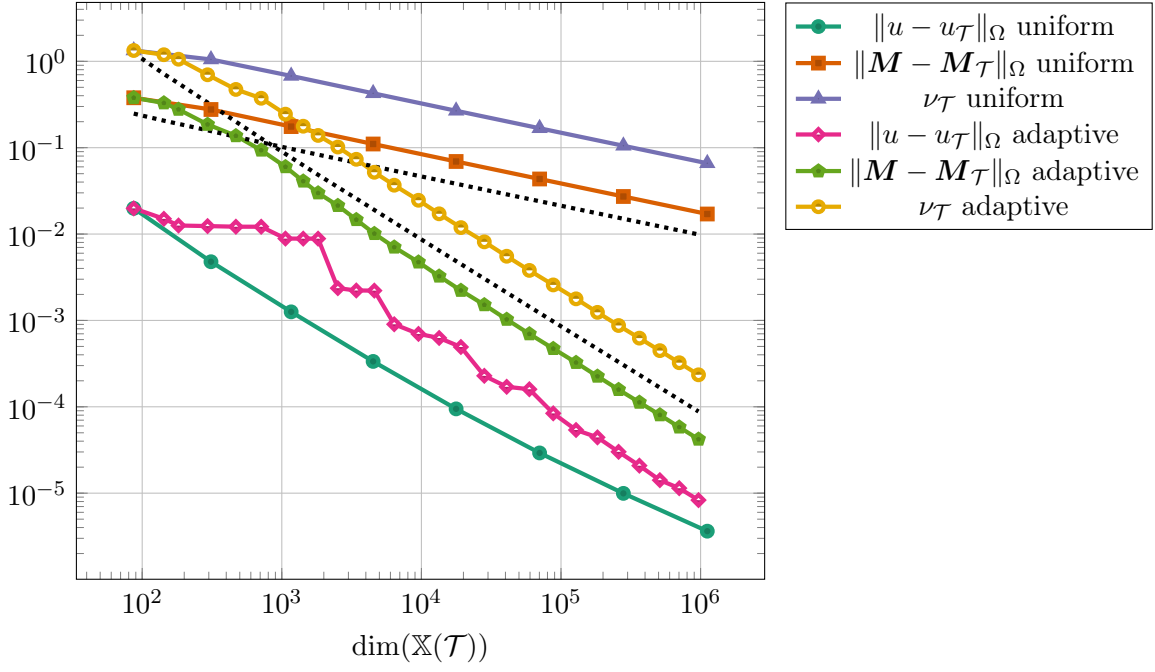
\begin{figure}
  \begin{center}
    \input{ExampleSingAdap}
    \caption{Errors for the singular problem from Section~\ref{sec:ex:sing}
              on sequences of uniformly and adaptively refined triangular meshes.
    The two black dotted lines indicate $\OO\big((\#\TT)^{-1}\big)$ and $\OO\big((\#\TT)^{-0.337}\big)$.}
    \label{fig:sing:adap}
  \end{center}
\end{figure}

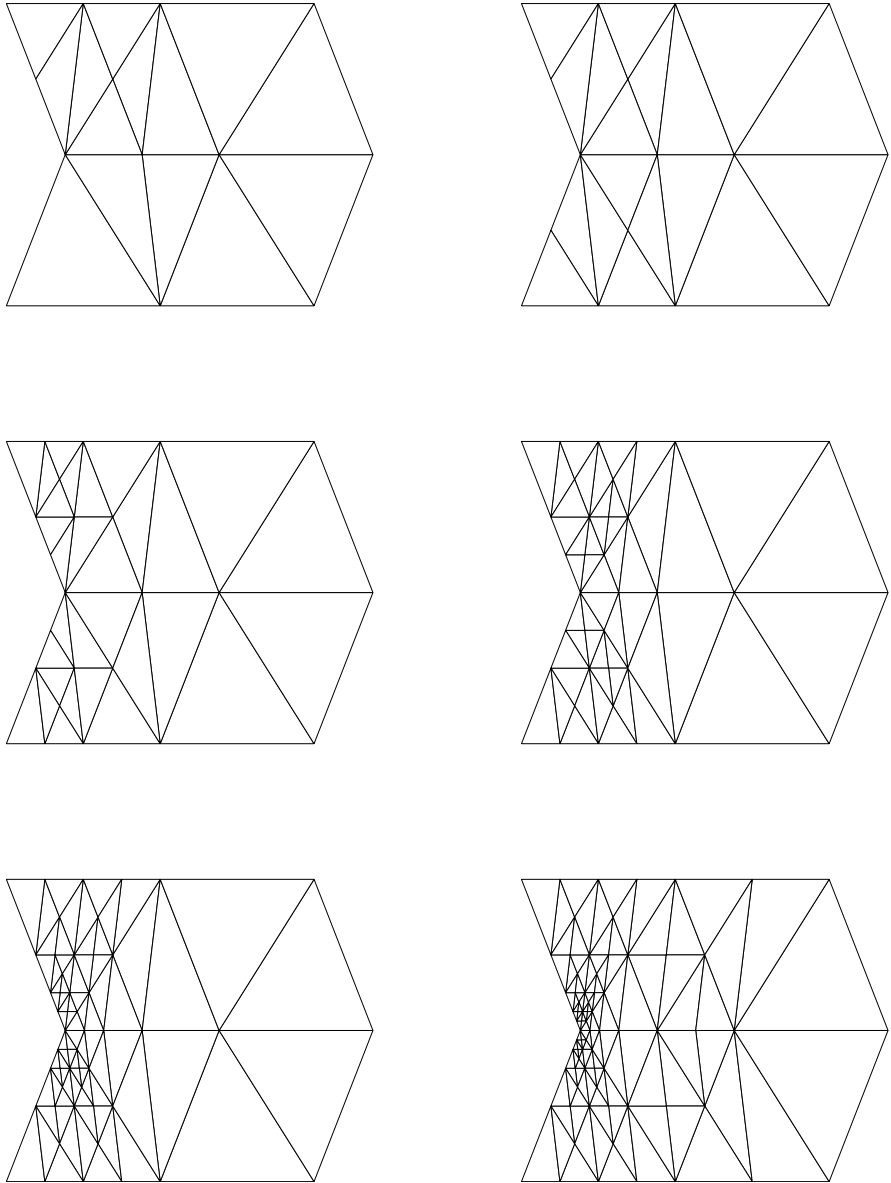
\begin{figure}
  \begin{center}
    \input{ExampleSingMeshes}
    \caption{Meshes $\TT_1,\ldots,\TT_6$ generated by the adaptive loop for the problem from Section~\ref{sec:ex:sing}. Axes are not scaled uniformly for presentation purposes.}
    \label{fig:sing:meshes}
  \end{center}
\end{figure}

Employing the adaptive loop described above with bulk criterion $\theta = 0.4$, we find that
convergence order $\OO\big( (\#\TT)^{-1}\big)$ is recovered, the one that we have seen before
for smooth solutions on uniformly refined meshes.
Figure~\ref{fig:sing:adap} shows the estimator and errors for uniformly and adaptively
refined triangulations. It illustrates the reliability and efficiency of the estimator in both cases.
Finally, Figure~\ref{fig:sing:meshes} shows the triangulations generated by the adaptive loop.
As expected, we observe strong refinements towards the re-entrant corner.

\subsection{Postprocessing}\label{sec:ex:postproc}
In this section we numerically study accuracy enhancement by the postprocessing scheme from Section~\ref{sec:postproc}. 
We consider the previously considered problems in convex and non-convex domains on triangular grids. 
Note that for the example from Section~\ref{sec:ex:smooth} the assumptions of Theorem~\ref{thm:postproc} are satisfied. The expected convergence order $\|u-u_\TT^\star\|_\Omega=\OO( (\#\TT)^{-2})$
is indeed observed in Figure~\ref{fig:postproc} (left plot). 
In contrast, the example from Section~\ref{sec:ex:sing} does not satisfy the assumptions of Theorem~\ref{thm:postproc}.
Therefore, improved convergence rates are not to be expected in general, see also Figure~\ref{fig:postproc} (right plot) for the case of a sequence of uniformly refined meshes. 
Interestingly, using adaptivity, optimal rates are recovered in the experiment, cf. Figure~\ref{fig:postproc}.

\begin{figure}
  \begin{center}
    \input{ExamplePostProc}
    \caption{Error of postprocessed solution for a problem with smooth solution (left) and a problem with singular solution (right). The two black dotted lines indicate $\OO\big((\#\TT)^{-1}\big)$ and $\OO\big((\#\TT)^{-2}\big)$.}
    \label{fig:postproc}
  \end{center}
\end{figure}

\appendix
\section{Basis functions}\label{sec:appendix}

\subsection{Basis functions for $\XXX_\triangle(K)$}\label{sec:appendix:triangle}
In this section we construct a basis for the local space $\XXX(K)$ on a triangle $K$.
Recall that $\Kref=\Kref_\triangle$ is the reference element with vertices $\widehat z_1=(0,0)$, $\widehat z_2=(1,0)$, $\widehat z_3 =(0,1)$, and edges $\widehat E_j$ spanned by $\widehat z_j$, $\widehat z_{j+1}$. We use cyclic indexing, i.e., $\widehat z_{j+3} = \widehat z_j$. 
For the physical element we use the analogous notation without $\widehat{(\cdot)}$.
We denote by $\widehat\tangential_j$ tangential vectors (positive orientation), and $\widehat\normal_j$ are the normal vectors on edge $\widehat E_j$.
Furthermore, $\widehat\eta_j$, $\eta_j$ denote the barycentric coordinate functions.

Defining {\footnotesize
\begin{alignat*}{2}
  \widehat\bT^\triangle_{(1,0)} &= \begin{pmatrix} -x^2y & -xy(y-1) \\ -xy(y-1) & -y(y-1)^2 \end{pmatrix},
  &\quad& \widehat\bT^\triangle_{(1,1)} = \begin{pmatrix} x^2(2x + y - 2) & x(2xy - 2y - x + y^2 + 1) \\ x(2xy - 2y - x + y^2 + 1) & y(y - 1)(2x + y - 1) \end{pmatrix}, \\
  \widehat\bT^\triangle_{(2,0)} &= \begin{pmatrix} x^2(x + y - 1) & xy(x + y - 1) \\ xy(x + y - 1)& y^2(x + y - 1) \end{pmatrix},
  &\quad& \widehat\bT^\triangle_{(2,1)} = \begin{pmatrix} x^2(y - x + 1) & -xy(x - y) \\ -xy(x - y) & -y^2(x - y + 1) \end{pmatrix}, \\
  \widehat\bT^\triangle_{(3,0)} &= \begin{pmatrix} -x(x - 1)^2 &  -xy(x - 1) \\ -xy(x - 1) & -xy^2 \end{pmatrix},
  &\quad& \widehat\bT^\triangle_{(3,1)} = \begin{pmatrix} -x(x - 1)(x + 2y - 1)& -y(2xy - y - 2x + x^2 + 1) \\ -y(2xy - y - 2x + x^2 + 1)  &  -y^2(x + 2y - 2) \end{pmatrix}
\end{alignat*}}
one verifies that
\begin{align*}
  \dual{\trdDiv{\Kref,\Eref_\ell,\tangential}\widehat\bT^\triangle_{(j,k)}}{\ell_{\widehat E_{\ell},m}}_{\widehat E_\ell} = \delta_{j,\ell}\delta_{k,m}
  \quad\text{and}\quad \widehat\bT^\triangle_{(j,k)}(\widehat z_\ell) = 0 = \trdDiv{\Kref,\Eref_\ell,\normal}\widehat\bT^\triangle_{(j,k)}.
\end{align*}
\begin{lemma}\label{lem:appendix:T}
  The transformed tensors $\bT^\triangle_{(j,k)} = \pk_K(\widehat\bT^\triangle_{(j,k)})\in\XXX(K)$ satisfy
  \begin{align*}
    \dual{\trdDiv{K,E_\ell,\tangential}\bT^\triangle_{(j,k)}}{\ell_{E_{\ell},m}}_{E_\ell} = \delta_{j,\ell}\delta_{k,m}
  \quad\text{and}\quad \bT^\triangle_{(j,k)}(z_\ell) = 0 = \trdDiv{K,E_\ell,\normal}\bT^\triangle_{(j,k)}.
  \end{align*}
\end{lemma}
\begin{proof}
  With $\widehat\bT^\triangle_{j,k}(\widehat z_\ell) = 0$ it follows that $\bT^\triangle_{j,k}(z_\ell)=0$ by definition of the transformation.

  Next, we show that $\trdDiv{K,E_\ell,\normal}\bT^\triangle_{(j,k)} = 0$. 
  Let $\ell=1$ and $p\in P^1(E)$.
  Take $v \in H^2(K)$ with $v|_{\partial K} = 0$, $\partial_{\normal} v|_{E_1} = \eta_{E_1} p$, $\partial_{\normal} v|_{E_m} = 0$ for $m=2,3$. Note that $\eta_{E_1} = \eta_1\eta_2$ is the edge bubble function.
  Then, 
  \begin{align*}
    -\dual{\trdDiv{K,E_1,\normal}\bT^\triangle_{(j,k)}}{\partial_\normal v}_{E_1} &= \dual{\trdDiv{K}\bT^\triangle_{(j,k)}}v_{\partial K} 
    = \dual{\trdDiv{\Kref}\widehat\bT^\triangle_{(j,k)}}{\widehat v}_{\partial\Kref} = 0.
  \end{align*}
  The latter identity follows because $\widehat v|_{\partial\Kref} = (v\circ F_K)|_{\partial\Kref} = 0$ and $\trdDiv{\Kref,\Eref_m,\normal}\widehat\bT^\triangle_{(j,k)} = 0$.
  Since $p\in P^1(E_1)$ and $\trdDiv{K,E_1,\normal}\bT^\triangle_{(j,k)}\in P^1(E)$ we conclude that $\trdDiv{K,E_1,\normal}\bT^\triangle_{(j,k)}=0$. For $\ell=2,3$ one argues similarly. 

  It remains to prove that $\dual{\trdDiv{K,E_\ell,\tangential}\bT^\triangle_{(j,k)}}{\ell_{E_{\ell},m}}_{E_\ell} = \delta_{j,\ell}\delta_{k,m}$.
  Let $\ell=1$. Let $p_m\in P^1(K)$ be such that $v_m=\eta_{E_1}p_m$ satisfies $\dual{\ell_{E_1,k}}{v_m}_{E_1} = \dual{\ell_{E_1,k}}{\ell_{E_1,m}}_{E_1}$ for all $k,m\in\{0,1\}$. 
  It follows that $\dual{q}{v_m}_{E_1} = \dual{q}{\ell_{E_1,m}}_{E_1}$ for all $q\in P^1(E_1)$ and $\dual{\widehat q}{\widehat v_m}_{\Eref_1} = \dual{\widehat q}{\ell_{\Eref_1,m}}_{\Eref_1}$ for all $\widehat q\in P^1(\Eref_1)$ with $\widehat v_m = v_m\circ F_K$.
  Then, 
  \begin{align*}
    \dual{\trdDiv{K,E_1,\tangential}\bT^\triangle_{(j,k)}}{\ell_{E_{1},m}}_{E_1} &= \dual{\trdDiv{K,E_1,\tangential}\bT^\triangle_{(j,k)}}{v_m}_{E_1}
    = \dual{\trdDiv{K}\bT^\triangle_{(j,k)}}{v_m}_{\partial K} \\
    &= \dual{\trdDiv{\Kref}\widehat\bT^\triangle_{(j,k)}}{\widehat v_m}_{\partial\Kref} 
    = \dual{\trdDiv{\Kref,\Eref_1,\tangential}\widehat\bT^\triangle_{(j,k)}}{\widehat v_m}_{\Eref_1}
    \\
    &= \dual{\trdDiv{\Kref,\Eref_1,\tangential}\widehat\bT^\triangle_{(j,k)}}{\ell_{\Eref_{1},m}}_{\Eref_1} = \delta_{j,1}\delta_{k,m}.
  \end{align*}
  For $\ell=2,3$ one argues similarly, which finishes the proof.
\end{proof}

For the basis functions associated to jump degrees of freedom we define
\begin{align*}
  \widetilde\bJ^\triangle_j = \frac2{\widehat\tangential_j\cdot\widehat\normal_{j+2}-\widehat\tangential_{j+2}\cdot\widehat\normal_j}\sym(\widehat\tangential_{j}\widehat\tangential_{j+2}^\top)\widehat\eta_{j}.
\end{align*}
One verifies that $\jump{\tangential\cdot\widetilde\bJ^\triangle_j\normal}_{\partial\Kref}(\widehat z_k) = \delta_{j,k}$ and $\trdDiv{\Kref,\Eref_k,\normal}\widetilde\bJ^\triangle_j = 0$.
However, the traces $\trdDiv{\Kref,\Eref_k,\tangential}\widetilde\bJ^\triangle_j$ do not vanish on all edges. To overcome this we subtract some correction terms: 
\begin{align*}
  \widehat\bJ^\triangle_j := \widetilde\bJ^\triangle_j - \sum_{\ell=1}^3 \dual{\trdDiv{\Kref,\Eref_\ell,\tangential}\widetilde\bJ^\triangle_j}{1}_{\Eref_\ell} \widehat\bT^\triangle_{\ell,0}.
\end{align*}

\begin{lemma}\label{lem:appendix:J}
  The transformed tensors $\bJ^\triangle_{j} = \pk_K(\widehat\bJ^\triangle_j)\in\XXX(K)$ satisfy
  \begin{align*}
    \trdDiv{K,E_\ell,\normal}\bJ^\triangle_j = 0 = \trdDiv{K,E_\ell,\tangential}\bJ^\triangle_j
    \quad\text{and}\quad 
    \jump{\tangential\cdot\bJ^\triangle_j\normal}_{\partial K} (z_k) = \delta_{j,k}.
  \end{align*}
\end{lemma}
\begin{proof}
  The identities $\trdDiv{K,E_\ell,\normal}\bJ^\triangle_j = 0 = \trdDiv{K,E_\ell,\tangential}\bJ^\triangle_j$ can be shown as in the proof of Lemma~\ref{lem:appendix:T}.
  For the last identity, take $v=\eta_k\in H^2(K)$. Then, 
  \begin{align*}
    -\delta_{j,k} &= -\jump{\tangential\cdot\widehat\bJ^\triangle_j\normal}_{\partial\Kref}(\widehat z_k)\widehat v(\widehat z_k) 
    = \dual{\trdDiv{\Kref}\widehat\bJ^\triangle}{\widehat v}_{\partial\Kref} 
    = \dual{\trdDiv{K}\bJ^\triangle}{v}_{\partial K} 
    \\ 
    &= -\jump{\tangential\cdot\bJ^\triangle_j\normal}_{\partial K}(z_k)v(z_k) = -\jump{\tangential\cdot\bJ^\triangle_j\normal}_{\partial K}(z_k),
  \end{align*}
  which finishes the proof.
\end{proof}

We define the remaining basis functions directly on $K$ (rather than by transformation from the reference element), 
\begin{align*}
  \widetilde\bN^\triangle_{(j,0)} = \frac{1}{\tangential_{j+1}\cdot\normal_j \, \tangential_{j+2}\cdot\normal_j}\sym(\tangential_{j+1}\tangential_{j+2}^\top), \quad
  \widetilde\bN^\triangle_{(j,1)} = \widetilde\bN^\triangle_{(j,0)}(\eta_{j+1}-\eta_j), \quad j=1,2,3.
\end{align*}
These functions satisfy $\trdDiv{K,E_\ell,\normal}\widetilde\bN^\triangle_{(j,k)} = \delta_{j,\ell}\ell_{E_j,k}$ but the other trace terms do not vanish in general. 
As before we subtract correction terms and define
\begin{align*}
  \bN^\triangle_{(j,0)} &= \widetilde\bN^\triangle_{(j,0)}-\sum_{\ell=1}^3 \jump{\tangential\cdot\widetilde\bN^\triangle_{(j,0)}\normal}_{\partial K}(z_\ell) \bJ^\triangle_\ell, \\
  \bN^\triangle_{(j,1)} &= \widetilde\bN^\triangle_{(j,1)}-\sum_{\ell=1}^3 \dual{\trdDiv{K,E_\ell,\tangential}\widetilde\bN^\triangle_{(j,1)}}{1}_{E_\ell} \bT^\triangle_{(\ell,0)} -\sum_{\ell=1}^3 \jump{\tangential\cdot\widetilde\bN^\triangle_{(j,1)}\normal}_{\partial K}(z_\ell) \bJ^\triangle_\ell.
\end{align*}
In the following lemma we collect some properties.
\begin{lemma}\label{lem:appendix:N}
  The tensors $\bN^\triangle_{(j,k)}\in\XXX(K)$ satisfy
  \begin{align*}
    \trdDiv{K,E_\ell,\tangential}\bN^\triangle_{(j,k)} = 0 = \jump{\tangential\cdot\bN^\triangle_j\normal}_{\partial K} (z_k)
    \quad\text{and}\quad 
    \trdDiv{K,E_\ell,\normal}\bN^\triangle_{(j,k)} = \delta_{j,\ell} \ell_{E_j,k}.
  \end{align*}
  Furthermore, 
  \begin{align*}
    \bN^\triangle_{(j,1)} = \widetilde\bN^\triangle_{(j,1)} - \sum_{\ell=1}^3 \jump{\tangential\cdot\widetilde\bN^\triangle_{(j,1)}\normal}_{\partial K}(z_\ell)\widetilde\bJ^{\triangle,K}_\ell \in \PPPsym^1(K)
\end{align*}
  where 
  \begin{align*}
    \widetilde\bJ_\ell^{\triangle,K} = \frac{2}{\tangential_j\cdot\normal_{j+2}-\tangential_{j+2}\cdot\normal_j} \sym(\tangential_j\tangential_{j+2}^\top)\eta_j.
\end{align*}
\end{lemma} 
\begin{proof}
  All but the last assertion follow by construction. 
  To see the last assertion we note that $\widetilde\bJ^{\triangle,K}_j$ is constructed as $\widetilde\bJ^\triangle_j$ but now on the element $K$ instead of $\Kref$. 
  As before one verifies that $\jump{\tangential\cdot \widetilde\bJ^{\triangle,K}_j\normal}_{\partial K}(z_k) = \delta_{j,k}$ and $\trdDiv{K,E_k,\normal}\widetilde\bJ^{\triangle,K}_j = 0$.
  Set 
  \begin{align*}
    \overline\bN^\triangle_{(j,1)} = \widetilde\bN^\triangle_{(j,1)} - \sum_{\ell=1}^3 \jump{\tangential\cdot\widetilde\bN^\triangle_{(j,1)}\normal}_{\partial K}(z_\ell)\widetilde\bJ^{\triangle,K}_\ell.
  \end{align*}
  The proof is finished if we show that $\overline\bN^\triangle_{(j,1)} = \bN^\triangle_{(j,1)}$. First, note that $\trdDiv{K,E_k,\normal} \overline\bN^\triangle_{(j,1)} = \trdDiv{K,E_k,\normal}\bN^\triangle_{(j,1)}$.
  Second, $\jump{\tangential\cdot\overline\bN^\triangle_{(j,1)}\normal}_{\partial K}(z_k) = 0$. Finally, $\trdDiv{K,E_k,\tangential}\overline\bN^\triangle_{(j,1)} = 0$ follows from a straightforward and simple but rather lengthy calculation (not shown).
\end{proof}

The following theorem summarizes the results of this section.
\begin{theorem}\label{thm:basisfun:triangle}
  The elements of $\mathbb{B}:=\{\bN^\triangle_{(j,k)},\bT^\triangle_{(j,k)}, \bJ^\triangle_j\,:\,j=1,2,3, \, k=0,1\}$ define a basis of $\XXX(K)$ and are dual to the degrees of freedom~\eqref{eq:dofs}.
  Moreover, 
  \begin{align*}
    \norm{\boldsymbol{B}}K \eqsim h_K \quad\forall \boldsymbol{B}\in \mathbb{B}.
  \end{align*}
\end{theorem}
\begin{proof}
  By construction the elements of $\mathbb{B}$ are dual to the degrees of freedom~\eqref{eq:dofs} and $\mathbb{B}\subset \XXX(K)$.
  We conclude that $\mathbb{B}$ is a basis of $\XXX(K)$.
  For the scaling properties we note that $\bT^\triangle_{(j,k)}$ and $\bJ^\triangle_j$ have been defined through the Piola--Kirchhoff transformation and by Lemma~\ref{lem:piolak} we get
  \begin{align*}
    \norm{\bT^\triangle_{(j,k)}}K \eqsim h_K\norm{\widehat\bT^\triangle_{(j,k)}}K \eqsim h_K \quad\text{and}\quad \norm{\bJ^\triangle_j}K \eqsim h_K\norm{\widehat\bJ^\triangle_j}K \eqsim h_K.
  \end{align*}
  To see that $\norm{\bN^\triangle_{(j,k)}}K\eqsim h_K$ set $\alpha_{(j,k),\ell} = \jump{\tangential\cdot\widetilde\bN_{(j,k)}^\triangle\normal}(z_\ell)$ and observe that $\sum_{\ell=1}^3 |\alpha_{(j,k),\ell}| \lesssim C$ with $C>0$ a generic constant independent of $h_K$, $j$ and $k$.
    Standard scaling arguments then show 
    \begin{align*}
      \norm{\bN_{(j,k)}^\triangle}K \eqsim h_K \norm{\pk_K^{-1}\bN_{(j,k)}^\triangle}{\Kref} 
      \eqsim h_K\big(1+\sum_{\ell=1}^3 |\alpha_{(j,k),\ell}|\big) \eqsim h_K,
    \end{align*}
    which finishes the proof.
\end{proof}

\subsection{Basis functions for $\XXX_\square(K)$}\label{sec:appendix:square}
In this section we construct a basis of the local space $\XXX(K)$ for a parallelogram $K$.
Recall that $\Kref=\Kref_\square$ is the reference element with vertices $\widehat z_1=(0,0)$, $\widehat z_2=(1,0)$, $\widehat z_3 =(1,1)$, $\widehat z_4=(0,1)$ and edges $\widehat E_j$ spanned by $\widehat z_j$, $\widehat z_{j+1}$. 
We adopt the notation from Section~\ref{sec:appendix:triangle} with obvious modifications. 
A slight difference is that here $\widehat\eta_j$ denotes the bilinear function on the reference element with $\eta_j(z_k) = \delta_{jk}$ and $\eta_j \circ F_K = \widehat\eta_j$.

Defining {\tiny
\begin{alignat*}{2}
  \widehat\bT^\square_{(1,0)} &= \begin{pmatrix} x(2y - 1)(x - 1)/4 & -y(2x - 1)(y - 1)/8 \\ -y(2x - 1)(y - 1)/8 & -y(y - 1)^2 \end{pmatrix},
  &\quad& \widehat\bT^\square_{(1,1)} = \begin{pmatrix} 0 & 3x(2y - 3)(x - 1)/8 \\ 3x(2y - 3)(x - 1)/8  & 3y(2x - 1)(y - 1)/4 \end{pmatrix}, \\
  \widehat\bT^\square_{(2,0)} &= \begin{pmatrix} x^2(x - 1) & x(2y - 1)(x - 1)/8 \\ x(2y - 1)(x - 1)/8 & -y(2x - 1)(y - 1)/4 \end{pmatrix},
  &\quad& \widehat\bT^\square_{(2,1)} = \begin{pmatrix} 3x(2y - 1)(x - 1)/4 & 3y(2x + 1)(y - 1)/8\\ 3y(2x + 1)(y - 1)/8 & 0 \end{pmatrix}, \\
  \widehat\bT^\square_{(3,0)} &= \begin{pmatrix} -x(2y - 1)(x - 1)/4 & y(2x - 1)(y - 1)/8  \\ y(2x - 1)(y - 1)/8 &  y^2(y - 1) \end{pmatrix},
  &\quad& \widehat\bT^\square_{(3,1)} = \begin{pmatrix} 0 & -3x(2y + 1)(x - 1)/8 \\ -3x(2y + 1)(x - 1)/8  &  -3y(2x - 1)(y - 1)/4 \end{pmatrix}, \\
  \widehat\bT^\square_{(4,0)} &= \begin{pmatrix} -x(x - 1)^2 & -x(2y - 1)(x - 1)/8  \\ -x(2y - 1)(x - 1)/8 & y(2x - 1)(y - 1)/4 \end{pmatrix},
  &\quad& \widehat\bT^\square_{(4,1)} = \begin{pmatrix} -3x(2y - 1)(x - 1)/4 &  -3y(2x - 3)(y - 1)/8 \\ -3y(2x - 3)(y - 1)/8  & 0  \end{pmatrix}
\end{alignat*}}
one verifies that
\begin{align*}
  \dual{\trdDiv{\Kref,\Eref_\ell,\tangential}\widehat\bT^\square_{(j,k)}}{\ell_{\widehat E_{\ell},m}}_{\widehat E_\ell} = \delta_{j,\ell}\delta_{k,m}
  \quad\text{and}\quad \widehat\bT^\square_{(j,k)}(\widehat z_\ell) = 0 = \trdDiv{\Kref,\Eref_\ell,\normal}\widehat\bT^\square_{(j,k)}.
\end{align*}
\begin{lemma}\label{lem:appendix:T:rect}
  The transformed tensors $\bT^\square_{(j,k)} = \pk_K(\widehat\bT^\square_{(j,k)})\in\XXX(K)$ satisfy
  \begin{align*}
    \dual{\trdDiv{K,E_\ell,\tangential}\bT^\square_{(j,k)}}{\ell_{E_{\ell},m}}_{E_\ell} = \delta_{j,\ell}\delta_{k,m}
  \quad\text{and}\quad \bT^\square_{(j,k)}(z_\ell) = 0 = \trdDiv{K,E_\ell,\normal}\bT^\square_{(j,k)}.
  \end{align*}
\end{lemma}
\begin{proof}
  The proof follows along the lines of the proof of Lemma~\ref{lem:appendix:T} and is therefore omitted.
\end{proof}

For the basis functions associated to jump degrees of freedom we define
\begin{align*}
  \widetilde\bJ^\square_j = \sym(\widehat\tangential_j\widehat\tangential_{j+1}^\top)\widehat\eta_j =  \frac{(-1)^j}2 \widehat\eta_j \begin{pmatrix} 0 & 1 \\ 1 & 0 \end{pmatrix}.
\end{align*}
One verifies that $\jump{\tangential\cdot\widetilde\bJ^\square_j\normal}_{\partial\Kref}(\widehat z_k) = \delta_{j,k}$ and $\trdDiv{\Kref,\Eref_k,\normal}\widetilde\bJ^\square_j = 0$. 
However, the traces $\trdDiv{\Kref,\Eref_k,\tangential} \widetilde\bJ^\square_j$ do not vanish on all edges. 
As in the previous section we subtract correction terms
\begin{align*}
  \widehat\bJ^\square_j := \widetilde\bJ^\square_j - \sum_{\ell=1}^4 \dual{\trdDiv{\Kref,\Eref_\ell,\tangential}\widetilde\bJ^\square_j}{1}_{\Eref_\ell}\widehat\bT^\square_{(\ell,0)} = \widetilde\bJ^\square_j - \widehat\bT^\square_{(j-1,0)} - \widehat\bT^\square_{(j,0)}. 
\end{align*}

\begin{lemma}\label{lem:appendix:J:rect}
  The transformed tensors $\bJ^\square_{j} = \pk_K(\widehat\bJ^\square_j)\in\XXX(K)$ satisfy
  \begin{align*}
    \trdDiv{K,E_\ell,\normal}\bJ^\square_j = 0 = \trdDiv{K,E_\ell,\tangential}\bJ^\square_j
    \quad\text{and}\quad 
    \jump{\tangential\cdot\bJ^\square_j\normal}_{\partial K} (z_k) = \delta_{j,k}.
  \end{align*}
\end{lemma}
\begin{proof}
  We only have to prove that $\bJ^\square_{j}\in\XXX(K)$. The proof of the other assertions follows along the lines of the proof of Lemma~\ref{lem:appendix:J} and is therefore omitted.
  Let $\mathbb{Q}^1(\Kref)\subset \PPPsym^2(\Kref)$ denote the space of symmetric tensor-valued polynomials with bilinear components. From the definition of $\XXX(\Kref)$ we find that $\mathbb{Q}^1(\Kref)\subset \XXX(\Kref)$.
Noting that $\widetilde\bJ_j^\square\in \mathbb{Q}^1(\Kref)$ we conclude that $\pk_K\widetilde\bJ_j^\square\in\XXX(K)$ by Proposition~\ref{propSpace:rect}, and, consequently, $\bJ_j^\square\in \XXX(K)$.
\end{proof}

As in the previous section, the remaining basis functions are defined directly on the physical element $K$, 
\begin{align*}
  \widetilde\bN^\square_{(j,0)} = \frac{1}{|\tangential_{j+1}\cdot\normal_j|^2} \tangential_{j+1}\tangential_{j+1}^\top(\eta_j+\eta_{j+1}), \quad
  \widetilde\bN^\square_{(j,1)} = \frac{1}{|\tangential_{j+1}\cdot\normal_j|^2} \tangential_{j+1}\tangential_{j+1}^\top(\eta_{j+1}-\eta_{j}).
\end{align*}
By construction these functions satisfy $\trdDiv{K,E_\ell,\normal}\widetilde\bN^\square_{(j,k)} = \delta_{j,\ell}\ell_{E_j,k}$.
As before we subtract correction terms to ensure that the other trace terms vanish, 
\begin{align*}
  \bN^\square_{(j,0)} &= \widetilde\bN^\square_{(j,0)} - \sum_{\ell=1}^4 \sum_{k=0}^1 \dual{\trdDiv{K,E,\tangential}\widetilde\bN^\square_{(j,0)}}{\ell_{E_\ell,k}}_{E_\ell}\bT^\square_{(\ell,k)} 
  -\sum_{\ell=1}^4 \jump{\tangential\cdot\widetilde\bN^\square_{(j,0)}\normal}_{\partial K}(z_\ell)\bJ^\square_\ell, \\
  \bN^\square_{(j,1)} &= \widetilde\bN^\square_{(j,1)} - \sum_{\ell=1}^4 \sum_{k=0}^1 \dual{\trdDiv{K,E,\tangential}\widetilde\bN^\square_{(j,1)}}{\ell_{E_\ell,k}}_{E_\ell}\bT^\square_{(\ell,k)}
  - \sum_{\ell=1}^4 \jump{\tangential\cdot\widetilde\bN^\square_{(j,1)}\normal}_{\partial K}(z_\ell)\bJ^\square_\ell.
\end{align*}
In the following lemma we collect some properties.
\begin{lemma}
  Tensors $\bN^\square_{(j,k)}\in\XXX(K)$ satisfy
  \begin{align*}
    \trdDiv{K,E_\ell,\tangential}\bN^\square_{(j,k)} = 0 = \jump{\tangential\cdot\bN^\square_j\normal}_{\partial K} (z_k)
    \quad\text{and}\quad 
    \trdDiv{K,E_\ell,\normal}\bN^\square_{(j,k)} = \delta_{j,\ell} \ell_{E_j,k}.
  \end{align*}
 \end{lemma}
 \begin{proof}
   We only have to show that $\bN^\square_{(j,k)}\in\XXX(K)$.
   The remaining assertions follow by definition. 
   Using the space of symmetric, bilinear tensors $\mathbb{Q}^1(\Kref)$ from before,
   we note that $\pk_K^{-1}\widetilde\bN^\square_{(j,k)}\in \mathbb{Q}^1(\Kref)$.
   It follows that $\widetilde\bN^\square_{(j,k)}\in\XXX(K)$ by Proposition~\ref{propSpace:rect}. 
     We conclude that $\bN^\square_{(j,k)}\in\XXX(K)$.
 \end{proof}

\begin{theorem}\label{thm:basisfun:rect}
  The elements of $\mathbb{B}:=\{\bN^\square_{(j,k)},\bT^\square_{(j,k)}, \bJ^\square_j\,:\,j=1,2,3,4, \, k=0,1\}$ define a basis of $\XXX(K)$ and are dual to the degrees of freedom~\eqref{eq:dofs}.
  Moreover, 
  \begin{align*}
    \norm{\boldsymbol{B}}K \eqsim h_K \quad\forall \boldsymbol{B}\in \mathbb{B}.
  \end{align*}
\end{theorem}
\begin{proof}
  The proof follows as for Theorem~\ref{thm:basisfun:triangle} and is therefore omitted.
\end{proof}

%===================================================================================================
\bibliographystyle{siam}
\bibliography{literature}
%===================================================================================================

\end{document}

%% file: header.tex
\newtheorem{theorem}{Theorem}

\newtheorem{lemma}[theorem]{Lemma}

\newtheorem{proposition}[theorem]{Proposition}
\newtheorem{remark}[theorem]{Remark}

\newcommand{\patch}{\omega}
\newcommand{\Patch}{\Omega}
\newcommand{\Kref}{\widehat{K}}

\newcommand{\Eref}{\widehat{E}}

\newcommand{\HHH}{\mathbb{H}}

\newcommand{\XXX}{\mathbb{X}}

\newcommand{\LLLts}{\mathbb{L}_{2,\mathrm{sym}}}
\newcommand{\PPP}{\mathbb{P}}
\newcommand{\PPPsym}{\mathbb{P}_\mathrm{sym}}

\newcommand{\conv}{\operatorname{conv}}

\newcommand\osc{\mathrm{osc}}

\DeclareMathOperator{\linhull}{span}

\newcommand{\projSZ}{J}
\newcommand{\projSZtilde}{\widetilde{J}}

\newcommand{\ip}[2]{(#1\hspace*{.5mm},#2)}
\newcommand{\dual}[2]{\langle#1\hspace*{.5mm},#2\rangle}

\newcommand{\norm}[3][]{#1\|#2#1\|_{#3}}

\newcommand{\curl}{\operatorname{curl}}
\newcommand{\rot}{\operatorname{rot}}
\newcommand{\Rot}{\operatorname{\mathbf{rot}}}
\newcommand{\Curl}{\operatorname{\mathbf{Curl}}}
\renewcommand\div{\operatorname{div}}
\def\Div{\operatorname{\mathbf{div}}}

\newcommand{\jump}[1]{[#1]}

\newcommand{\Hdivset}[1]{\boldsymbol{H}(\div;#1)}
\newcommand{\HDivset}[1]{\mathbb{H}(\Div;#1)}

\newcommand{\HdivDivset}[1]{\mathbb{H}(\dDiv;#1)}

\newcommand{\piola}{\mathcal{P}}
\newcommand{\pk}{\mathcal{H}}
\newcommand{\Ctens}{\mathcal{C}}

\newcommand\trdDiv[1]{\operatorname{tr}_{#1}^{\mathrm{dDiv}}}
\newcommand\trGgrad[1]{\operatorname{tr}_{#1}^{\mathrm{Ggrad}}}

\newcommand{\set}[2]{\big\{#1\,:\,#2\big\}}

\newcommand{\RT}{\ensuremath{\boldsymbol{RT}}}

\newcommand{\R}{\ensuremath{\mathbb{R}}}
\newcommand{\N}{\ensuremath{\mathbb{N}}}
\newcommand{\HH}{\ensuremath{{\boldsymbol{H}}}}

\newcommand{\LL}{\ensuremath{\boldsymbol{L}}}

\newcommand{\TT}{\ensuremath{\mathcal{T}}}
\newcommand{\marked}{\ensuremath{\mathcal{M}}}

\newcommand{\OO}{\ensuremath{\mathcal{O}}}
\newcommand{\edges}{\ensuremath{\mathcal{E}}}
\newcommand{\vertices}{\ensuremath{\mathcal{V}}}

\newcommand{\normal}{\ensuremath{{\boldsymbol{n}}}}
\newcommand{\tangential}{\ensuremath{{\boldsymbol{t}}}}

\newcommand\dDiv{\operatorname{div\mathbf{div}}}

\newcommand\Grad{\boldsymbol\nabla}
\newcommand\grad{\nabla}
\newcommand\Ggrad{\Grad\grad}

\newcommand{\bphi}{{\boldsymbol\phi}}
\newcommand{\bpsi}{{\boldsymbol\psi}}

\newcommand{\bsigma}{{\boldsymbol\sigma}}

\newcommand{\bq}{{\boldsymbol{q}}}
\newcommand{\bu}{\boldsymbol{u}}

\newcommand{\sym}{\operatorname{sym}}
\newcommand{\symCurl}{\operatorname{\mathbf{sCurl}}}

\newcommand{\projHdDiv}{\Pi^{\dDiv}}

\newcommand{\ba}{\boldsymbol{a}}
\newcommand{\bb}{\boldsymbol{b}}
\newcommand{\bc}{\boldsymbol{c}}
\newcommand{\bd}{\boldsymbol{d}}
\newcommand{\bx}{\boldsymbol{x}}

\newcommand{\bg}{\boldsymbol{g}}

\newcommand{\bM}{\boldsymbol{M}}
\newcommand{\bI}{\boldsymbol{I}}
\newcommand{\bJ}{\boldsymbol{J}}
\newcommand{\bH}{\boldsymbol{H}}
\newcommand{\bP}{\boldsymbol{P}}

\newcommand{\bB}{\boldsymbol{B}}
\newcommand{\bN}{\boldsymbol{N}}
\newcommand{\bT}{\boldsymbol{T}}
\newcommand{\bQ}{\boldsymbol{Q}}
\newcommand{\bX}{\boldsymbol{X}}
\newcommand{\bY}{\boldsymbol{Y}}

%% file: ExampleSmooth.tex
\begin{tikzpicture}
\begin{loglogaxis}[
width=0.49\textwidth,
%cycle list name=black white,
%cycle list name=exotic,
cycle list/Dark2-6,
% combine it with ’mark list*’:
cycle multiindex* list={
mark list*\nextlist
Dark2-6\nextlist},
every axis plot/.append style={ultra thick},
xlabel={$\dim(\XXX(\TT))$},
grid=major,
legend entries={\tiny $\|u-u_\TT\|_\Omega$,\tiny $\|\bM-\bM_\TT\|_\Omega$,\tiny $\|\dDiv(\bM-\bM_\TT)\|_\Omega$},
legend pos=south west,
%legend pos=outer north east,
]
\addplot table [x=dof,y=errU] {data/ExampleSmooth.dat};
\addplot table [x=dof,y=errM] {data/ExampleSmooth.dat};
\addplot table [x=dof,y=errDDivM] {data/ExampleSmooth.dat};
\addplot [black,dotted,mark=none] table [x=dof,y expr={0.5*\thisrowno{0}^(-1)}] {data/ExampleSmooth.dat};
\end{loglogaxis}
\end{tikzpicture}
\begin{tikzpicture}
\begin{loglogaxis}[
width=0.49\textwidth,
%cycle list name=black white,
%cycle list name=exotic,
cycle list/Dark2-6,
% combine it with ’mark list*’:
cycle multiindex* list={
mark list*\nextlist
Dark2-6\nextlist},
every axis plot/.append style={ultra thick},
xlabel={$\dim(\XXX(\TT))$},
grid=major,
legend entries={\tiny $\|u-u_\TT\|_\Omega$,\tiny $\|\bM-\bM_\TT\|_\Omega$,\tiny $\|\dDiv(\bM-\bM_\TT)\|_\Omega$},
legend pos=south west,
% legend style={at={(0.5,1.05)},anchor=south},
% legend columns=3, 
%         legend style={
%                     % the /tikz/ prefix is necessary here...
%                     % otherwise, it might end-up with `/pgfplots/column 2`
%                     % which is not what we want. compare pgfmanual.pdf
%             /tikz/column 2/.style={
%                 column sep=5pt,
%             }},
%legend pos=outer north east,
]
\addplot table [x=dof,y=errU] {data/ExampleSmoothRect.dat};
\addplot table [x=dof,y=errM] {data/ExampleSmoothRect.dat};
\addplot table [x=dof,y=errDDivM] {data/ExampleSmoothRect.dat};
\addplot [black,dotted,mark=none] table [x=dof,y expr={0.5*\thisrowno{0}^(-1)}] {data/ExampleSmoothRect.dat};
\end{loglogaxis}
\end{tikzpicture}

%% file: ExampleSingInitMesh.tex
\begin{tikzpicture}
\begin{axis}[
width=0.49\textwidth,
    axis equal,
    xlabel = {$x$},
%    width=0.49\textwidth,
%    xlabel={$x$},
%    ylabel={$y$},
]%, axis equal image = true]

\addplot[patch,color=white,%mesh, 
faceted color = black, line width = 0.5pt,
patch table ={data/ele1tri.dat}] file{data/coo1tri.dat};
\end{axis}
\end{tikzpicture}
\begin{tikzpicture}
\begin{axis}[
width=0.49\textwidth,
    axis equal,
    xlabel = {$x$},
%    width=0.49\textwidth,
%    xlabel={$x$},
%    ylabel={$y$},
]%, axis equal image = true]

\addplot[patch,patch type = rectangle,color=white,%mesh, 
faceted color = black, line width = 0.5pt,
patch table ={data/ele1rect.dat}] file{data/coo1rect.dat};
\end{axis}
\end{tikzpicture}

%% file: ExampleSing.tex
\begin{tikzpicture}
\begin{loglogaxis}[
width=0.65\textwidth,
%cycle list name=black white,
%cycle list name=exotic,
cycle list/Dark2-6,
% combine it with ’mark list*’:
cycle multiindex* list={
mark list*\nextlist
Dark2-6\nextlist},
every axis plot/.append style={ultra thick},
xlabel={$\dim(\XXX(\TT))$},
grid=major,
legend entries={$\|u-u_\TT\|_\Omega$ (tri.),$\|\bM-\bM_\TT\|_\Omega$ (tri.),$\|u-u_\TT\|_\Omega$ (par.),$\|\bM-\bM_\TT\|_\Omega$ (par.)},
%legend pos=south west,
legend pos=outer north east,
]

\addplot table [x=dof,y=errU] {data/ExampleSingUnif.dat};
\addplot table [x=dof,y=errM] {data/ExampleSingUnif.dat};
\addplot table [x=dof,y=errU] {data/ExampleSingRect.dat};
\addplot table [x=dof,y=errM] {data/ExampleSingRect.dat};

\addplot [black,dotted,mark=none] table [x=dof,y expr={0.5*\thisrowno{0}^(-0.337)}] {data/ExampleSingUnif.dat};
\addplot [black,dotted,mark=none] table [x=dof,y expr={0.2*\thisrowno{0}^(-1)}] {data/ExampleSingUnif.dat};
\end{loglogaxis}
\end{tikzpicture}

%% file: ExampleSingAdap.tex
\begin{tikzpicture}
\begin{loglogaxis}[
width=0.65\textwidth,
%cycle list name=black white,
%cycle list name=exotic,
cycle list/Dark2-6,
% combine it with ’mark list*’:
cycle multiindex* list={
mark list*\nextlist
Dark2-6\nextlist},
every axis plot/.append style={ultra thick},
xlabel={$\dim(\XXX(\TT))$},
grid=major,
legend entries={$\|u-u_\TT\|_\Omega$ uniform,$\|\bM-\bM_\TT\|_\Omega$ uniform,$\nu_\TT$ uniform,$\|u-u_\TT\|_\Omega$ adaptive,$\|\bM-\bM_\TT\|_\Omega$ adaptive,$\nu_\TT$ adaptive},
%legend pos=south west,
legend pos=outer north east,
]
\addplot table [x=dof,y=errU] {data/ExampleSingUnif.dat};
\addplot table [x=dof,y=errM] {data/ExampleSingUnif.dat};
\addplot table [x=dof,y=estM] {data/ExampleSingUnif.dat};
\addplot table [x=dof,y=errU] {data/ExampleSingAdap.dat};
\addplot table [x=dof,y=errM] {data/ExampleSingAdap.dat};
\addplot table [x=dof,y=estM] {data/ExampleSingAdap.dat};

\addplot [black,dotted,mark=none] table [x=dof,y expr={0.5*\thisrowno{0}^(-0.337)}] {data/ExampleSingAdap.dat};
\addplot [black,dotted,mark=none] table [x=dof,y expr={10*\thisrowno{0}^(-1)}] {data/ExampleSingAdap.dat};
\end{loglogaxis}
\end{tikzpicture}

%% file: ExampleSingMeshes.tex
% \begin{tikzpicture}
% \begin{axis}[hide axis,
% %    title={$\#\TT=114$},
% width=0.5\textwidth,
%     axis equal,
% %    width=0.49\textwidth,
% %    xlabel={$x$},
% %    ylabel={$y$},
% ]%, axis equal image = true]
% 
%     \nextgroupplot
%         \addplot3[patch,shader=interp] table{data/HyperbolicSolExact.dat};
% 
% 
% \end{axis}
% \end{tikzpicture}
% 
\begin{tikzpicture}
  \begin{groupplot}[group style = {group size = 2 by 3},width=0.45\textwidth,hide axis]
%     \nextgroupplot
%         \addplot[patch,color=white,
%           faceted color = black, line width = 0.3pt,
%           patch table ={data/ele1.dat}] file{data/coo1.dat};
    \nextgroupplot
        \addplot[patch,color=white,%mesh,
          faceted color = black, line width = 0.3pt,
          patch table ={data/ele2.dat}] file{data/coo2.dat};
    \nextgroupplot
        \addplot[patch,color=white,%mesh,
          faceted color = black, line width = 0.3pt,
          patch table ={data/ele3.dat}] file{data/coo3.dat};
    \nextgroupplot
        \addplot[patch,color=white,%mesh,
          faceted color = black, line width = 0.3pt,
          patch table ={data/ele4.dat}] file{data/coo4.dat};
    \nextgroupplot
        \addplot[patch,color=white,%mesh,
          faceted color = black, line width = 0.3pt,
          patch table ={data/ele5.dat}] file{data/coo5.dat};
    \nextgroupplot
        \addplot[patch,color=white,%mesh,
          faceted color = black, line width = 0.3pt,
          patch table ={data/ele6.dat}] file{data/coo6.dat};
   \nextgroupplot
       \addplot[patch,color=white,%mesh,
         faceted color = black, line width = 0.3pt,
         patch table ={data/ele7.dat}] file{data/coo7.dat};
%    \nextgroupplot
%        \addplot[patch,color=white,%mesh,
%          faceted color = black, line width = 0.3pt,
%          patch table ={data/ele8.dat}] file{data/coo8.dat};
%    \nextgroupplot
%        \addplot[patch,color=white,%mesh,
%          faceted color = black, line width = 0.3pt,
%          patch table ={data/ele9.dat}] file{data/coo9.dat};
  \end{groupplot}
\end{tikzpicture}

%% file: ExamplePostProc.tex
\begin{tikzpicture}
\begin{loglogaxis}[
width=0.49\textwidth,
%cycle list name=black white,
%cycle list name=exotic,
cycle list/Dark2-6,
% combine it with ’mark list*’:
cycle multiindex* list={
mark list*\nextlist
Dark2-6\nextlist},
every axis plot/.append style={ultra thick},
xlabel={$\dim(\XXX(\TT))$},
grid=major,
legend entries={\tiny $\|u-u_\TT\|_\Omega$,\tiny $\|u-u_\TT^\star\|_\Omega$},
legend pos=south west,
%legend pos=outer north east,
]
\addplot table [x=dof,y=errU] {data/ExampleSmoothPostProc.dat};
\addplot table [x=dof,y=errUstar] {data/ExampleSmoothPostProc.dat};
\addplot [black,dotted,mark=none] table [x=dof,y expr={5e-2*\thisrowno{0}^(-1)}] {data/ExampleSmoothPostProc.dat};
\addplot [black,dotted,mark=none] table [x=dof,y expr={5e-2*\thisrowno{0}^(-2)}] {data/ExampleSmoothPostProc.dat};
\end{loglogaxis}
\end{tikzpicture}
\begin{tikzpicture}
\begin{loglogaxis}[
width=0.49\textwidth,
%cycle list name=black white,
%cycle list name=exotic,
cycle list/Dark2-6,
% combine it with ’mark list*’:
cycle multiindex* list={
mark list*\nextlist
Dark2-6\nextlist},
every axis plot/.append style={ultra thick},
xlabel={$\dim(\XXX(\TT))$},
grid=major,
legend entries={\tiny $\|u-u_\TT\|_\Omega$ unif.,\tiny $\|u-u_\TT^\star\|_\Omega$ unif.,\tiny $\|u-u_\TT\|_\Omega$ adap.,\tiny $\|u-u_\TT^\star\|_\Omega$ adap.},
legend pos=south west,
%legend pos=outer north east,
]
\addplot table [x=dof,y=errU] {data/ExampleSingularUnifPostProc.dat};
\addplot table [x=dof,y=errUstar] {data/ExampleSingularUnifPostProc.dat};
\addplot table [x=dof,y=errU] {data/ExampleSingularAdapPostProc.dat};
\addplot table [x=dof,y=errUstar] {data/ExampleSingularAdapPostProc.dat};
\addplot [black,dotted,mark=none] table [x=dof,y expr={6e0*\thisrowno{0}^(-1)}] {data/ExampleSingularAdapPostProc.dat};
\addplot [black,dotted,mark=none] table [x=dof,y expr={4e-1*\thisrowno{0}^(-2)}] {data/ExampleSingularAdapPostProc.dat};
\end{loglogaxis}
\end{tikzpicture}